\documentclass[11pt]{article}
 \pdfoutput=1

\usepackage{amsmath,amssymb,amsfonts,latexsym,amsbsy,epsfig, mathrsfs,
multirow, ifthen,extarrows}
\usepackage{ifthen}

\usepackage{algorithm}
\usepackage{algpseudocode}
\usepackage{amsmath,amssymb,amsfonts,latexsym,amsbsy,epsfig, mathrsfs,
multirow, ifthen,extarrows}
\usepackage{ifthen}
\usepackage{comment}
\usepackage{pgf,tikz}
\usepackage{amsmath,amssymb,mathtools}
\usepackage{stmaryrd}
\usepackage[mathscr]{eucal}
\usepackage{multirow}
\usepackage{enumerate}
\usepackage{subfig}
\usepackage{xr}

\usetikzlibrary{arrows, decorations.markings,positioning,swigs}

\usetikzlibrary{shapes}

\usepackage[authoryear]{natbib}
\usepackage{url}
\usepackage{hyperref}

\DeclareMathOperator{\pa}{pa}

\DeclareMathOperator{\an}{an}

\DeclareMathOperator{\pre}{pre}

\newtheorem{theorem}{Theorem}

\newtheorem{lemma}[theorem]{Lemma}

\newtheorem{claim}[theorem]{Claim}

\newtheorem{definition}[theorem]{Definition}

\newcommand{\G}{{\mathscr G}}

\newcommand\indep{\protect\mathpalette{\protect\independenT}{\perp}}
\def\independenT#1#2{\mathrel{\rlap{$#1#2$}\mkern2.5mu{#1#2}}}
\newcommand\indepd{\protect\mathpalette{\protect\independenTd}{\perp}}
\def\independenTd#1#2{\mathrel{\rlap{$#1#2$}\mkern2.5mu{#1#2}_d}}

\begin{document}

\title{
Potential Outcome and Decision Theoretic Foundations for Statistical Causality}
\author{Thomas S. Richardson\\ University of Washington \and James M. Robins\\ Harvard School of Public Health} 




\maketitle

\section{Introduction}

In his recent paper, {\it Decision Theoretic Foundations for Causality}', Philip Dawid elaborates on an earlier theory that he advanced previously \citep{dawid02influence}. We welcome Dawid's efforts to build a foundation for causal models that aims  to develop a graphical framework, while at the same time placing an emphasis on making assumptions that are both transparent and testable. Similar concerns have also motivated much of our previous work on potential outcome models represented in terms of Finest Fully Randomized Causally Interpreted Structured Tree Graphs \citep{robins86new}  and Single World Intervention Graphs
 \citep{thomas13swig}.

Indeed, like Dawid, we have argued that, in contrast, the assumption of independent errors that is typically adopted by users of Pearl's Non-Parametric Structural Equations (also called Structural Causal Models) is untestable and also imposes (super exponentially) many assumptions that are unnecessary for most purposes; further, the independent error assumptions allows the identification of causal quantities that cannot be identified via any randomized experiment on the observed variables \citep{robins10alternative}. Thus this assumption contradicts the dictum ``no causation without manipulation'' and severs the connection between experimentation and causal inference that has been central to much of the conceptual progress during the last century.
We also note that \citet{imbens:2022} cites the move to specifying causal models using 
potential outcomes rather than error terms as underpinning the `credibility revolution' in Econometrics.

In our view Dawid's updated theory represents a marked advance on his earlier proposal in that it requires stronger ontological commitments, specifically, the existence of an `intent to treat' (ITT) variable, before a model may be called causal. ITT variables are necessary and important in order to encode the notion of ignorability and the effect of treatment on the treated.

In addition, as noted by Dawid, the ITT variables make it possible to connect his approach to that based on potential outcomes\footnote{Though Dawid and others distinguish between potential outcomes and counterfactuals on philosophical grounds we do not do so here; we think that this distinction, though of interest, is a separate issue from those under discussion here.} and Single World Interventions Graphs. The connection between the two approaches may help to illuminate the strengths and weakness of each formalism.
We also present a reformulation of Dawid's theory that is essentially 
equivalent to his proposal and isomorphic to SWIGs.

We thank Philip Dawid for helpful feedback on our paper; in particular, for pointing out a significant omission regarding our proposed definition of distributional consistency for SWIGs. We also thank him for his patience regarding the completion of this manuscript.

\section{Relating Observational and Experimental Worlds}

At a high-level, every approach to causal inference relates a model describing a factual passively observed world and models describing hypothetical `interventional' worlds in which a treatment (or exposure) variable takes on a specific value. 

In both the current and previous decision-theoretic conceptions advocated by Dawid these worlds `exist' at least hypothetically, as different distributions. The relation is then created by the assertion of equalities linking different parts of these distributions. In Dawid's formalism the  set of distributions is represented using a single {\it kernel} object in which non-random regime indicators (also called `policy variables' by \citet{sgs:1993}) index the different distributions; there is no requirement that these distributions live on the same probability space.
Dawid encodes the equalities between the observational and interventional worlds via extended conditional independence relations, including independence from (and conditional on) regime indicators.

In the standard presentation of the potential outcome approach, random variables corresponding to the outcomes for an individual under all possible interventions\footnote{This does not mean that it is assumed that all variables can be intervened on.} are assumed to exist, living on a common probability space. 
The consistency assumption then serves to  construct the factual variables as a deterministic function of the potential outcomes.
Owing to the fundamental problem of causal inference the resulting factual distribution is consistent with many different intervention distributions.
However,  under additional Markov restrictions on the joint distribution of the potential outcomes,  the interventional distribution becomes identified from the joint distribution of the factuals under a positivity assumption.
Notwithstanding this, often in practice data are obtained on a subset of the factual variables in which case some or even all interventional distributions become only partially identified from the available (i.e.~the observed) data.

\begin{figure}
\centering
\begin{tikzpicture}
\tikzset{line width=1.25pt, outer sep=0pt,
         ell/.style={draw, inner sep=2pt,
          line width=1.5pt}, >=stealth,
         swig vsplit={gap=5pt,
         inner line width right=0pt, line color right=red},
         ell2/.style={draw,fill opacity=0, text opacity=1}};
\begin{scope}
\node[name=a,ell,draw,shape=ellipse]{$A$};
\node[name=b,ell,draw,shape=ellipse,right=12mm of a]{$B$};
\node[name=c,ell,draw,shape=ellipse,right=12mm of b]{$C$};
\draw[ell,->](a) to (b);
\draw[ell,->](b) to (c);
\node[name=atag, left=3mm of a]{(a)};
\end{scope}
\begin{scope}[yshift=-2cm,xshift=0]  
\node[name=btag, below= of atag]{(b)};
\node[name=a,right=3mm of btag, ell,shape=swig vsplit, ell2]{
        \nodepart{left}{$A(a,b)$}
        \nodepart{right}{$a$}};
\node[name=b,ell,right=8mm of a, shape=swig vsplit, ell2]{
        \nodepart{left}{$B(a,b)$}
        \nodepart{right}{$b$}};
\node[name=c,ell,draw,shape=ellipse,right=8mm of b]{$C(a,b)$};
\draw[ell,->](a) to (b);
\draw[ell,->](b) to (c);
\end{scope}
\begin{scope}[yshift=-4cm,xshift=0]  
\node[name=ctag, below= of btag]{(c)};
\node[name=a,right=3mm of ctag, ell,shape=swig vsplit, inner sep=3pt, ell2]{
        \nodepart{left}{$A$}
        \nodepart{right}{$a$}};
\node[name=b,ell,right=8mm of a, shape=swig vsplit, ell2]{
        \nodepart{left}{$B(a)$}
        \nodepart{right}{$b$}};
\node[name=c,ell,draw,shape=ellipse,right=8mm of b]{$C(a,b)$};
\draw[ell,->](a) to (b);
\draw[ell,->](b) to (c);
\end{scope}
\begin{scope}[yshift=-4cm,xshift=0]  
\node[name=dtag, below= of ctag]{(d)};
\node[name=a,right=3mm of dtag, ell,shape=swig vsplit, inner sep=3pt, ell2]{
        \nodepart{left}{$A$}
        \nodepart{right}{$a$}};
\node[name=b,ell,right=8mm of a, shape=swig vsplit, ell2]{
        \nodepart{left}{$B(a)$}
        \nodepart{right}{$b$}};
\node[name=c,ell,draw,shape=ellipse,right=8mm of b]{$C(b)$};
\draw[ell,->](a) to (b);
\draw[ell,->](b) to (c);
\end{scope}
\end{tikzpicture}
\caption{Illustration of SWIG labeling schemes.
(a) DAG $\G$ representing the observed joint distribution $p(A,B,C)$; (b) SWIG ${\G}(a,b)$ with uniform labeling;
(c) SWIG ${\G}(a,b)$ with temporal labeling; (d) SWIG ${\G}(a,b)$ with ancestral labeling.
(These and other figures were created using the \href{https://ctan.org/pkg/tikz-swigs?lang=en}{\tt swigs} Ti{\it k}Z package, available on \href{https://ctan.org}{CTAN}.) 
\label{fig:swig-label}}
\end{figure}
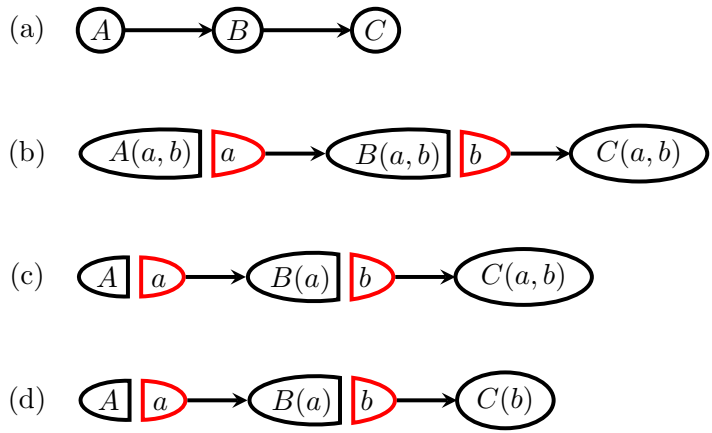

\section{Single-World Intervention Graphs}\label{sec:swigs}

The Single-World Intervention Graph (SWIG) approach is designed to provide a simple way to relate graphs representing joint distributions over the observed variables and those representing joint distributions over potential outcomes. The approach is `Single World' in that each of the constraints defining the model concerns a set of potential outcomes corresponding to a single joint intervention on the target variables.\footnote{Consequently, though SWIGs define a potential outcome model, they do not impose `cross-world' assumptions such as strong ignorability $Y(0),Y(1) \indep T$ or strong conditional ignorability, $Y(0),Y(1) \indep T \mid X$.}

Following \citep{robins86new,thomas13swig} we will assume throughout that there is a set of variables indexed by $V=\{1,\ldots , p\}$ and that a pre-specified (possibly strict)  subset $A \subseteq V$ of these variables are targets for intervention. Often we will, with slight abuse of notation, also refer to the corresponding sets of random variables as $V$ and $A$ respectively.

However, for proofs and formal statements it is sometimes necessary to distinguish between the random variables and the sets that index them. For this purpose we introduce the following notation: we define $X_B \equiv \{X_i, i \in B \subseteq V\}$, so that the complete set of factual variables is $X_V$ and the subset that are targets for intervention are $X_A$. We use $\mathfrak{X}_{i}$ as the state space for the variable $X_i$ and we will let  $\mathfrak{X}_V \equiv \times_{i \in V} \mathfrak{X}_{i}$  and $\mathfrak{X}_A \equiv \times_{i \in A} \mathfrak{X}_{i}$ be the state spaces for the variables with indices in $V$ and $A$ respectively. Similarly, given an assignment $x_V$ to the variables (with indices) in $V$, we let $x_i$ and $x_B$ refer to the value assigned to $X_i$ and to the set $X_B$. We also make use of the usual shorthand, using, for example $A_i$ to refer to $X_{A_i}$, $A$ for $X_A$  and $a_i$ to denote $x_{a_i}$.

\begin{definition}
Given a Directed Acyclic Graph (DAG) $\G$ with vertex set $V$, the SWIG $\G({a})$ corresponding to an intervention that sets the variables in $A=\{A_1,\ldots ,A_k\}\subseteq V$
to ${a}=(a_1,\ldots ,a_k) \in \mathfrak{X}_{A}$ is constructed as follows:
\begin{itemize}
\item[\rm (1)] Every vertex $A_i \in A$ is split into two halves, a `random half' and a `fixed half';
\item[\rm (2)] The random half contains $A_i$ and inherits all of the incoming edges directed into $A_i$ in the original graph;
\item[\rm (3)] The fixed half inherits all of the outgoing edges directed out of $A_i$ in the original graph and is labeled with the value $a_i$;
\item[\rm (4)] Random vertices in nodes on the graph are then re-labeled according to one of the schemes below.
\end{itemize}
\end{definition}

There are three labeling schemes that may be employed in step (4):
\begin{itemize}
\item[] {\bf Uniform labeling:} Every random vertex $Y$ in the SWIG $\G({a})$ is labeled with the full vector $Y({a_1,\ldots, a_k})$;
\item[] {\bf Temporal labeling:} Given a total ordering of the vertices on the original graph, each random vertex $Y$
is labeled $Y(a_1,\ldots a_i)$ with the values  corresponding to those vertices $A_1,\ldots , A_i$ that are ordered prior to $Y$;
\item[]{\bf Ancestral labeling:} Each random vertex $Y$ is labeled $Y(a_{\an_{\G(a)}(Y)})$ where ${a_{\an_{\G(a)}(Y)}}$ corresponds to those fixed vertices $a_i$ that are {\em still} ancestors of $Y$ {\em after} splitting the nodes in $A$.
\end{itemize}

Temporal labeling may be seen as encoding the assumption that interventions in the future do not affect outcomes in the past. Thus the potential outcome for a variable 
$Y(a_1,\ldots ,a_k)$, in a world in which there is an  intervention on $A_1,\ldots, A_k$,
is a function only of those interventions $A_1,\ldots, A_i$ that took place (temporally) before $Y$, so that $Y(a_1,\ldots ,a_k)= Y(a_1,\ldots ,a_i)$.
This is the natural labeling scheme to apply in the context where all variables are temporally ordered and missing edges correspond (solely) to the 
absence of population level direct effects.

Ancestral labeling encodes the assumption that the potential outcome $Y(a_1,\ldots ,a_k)$ is solely a function of those interventions that are (still) causally antecedent to $Y$ in the context of the other interventions that are being carried out. Thus, for example, in Figure \ref{fig:swig-label} (d), the vertex for $C$ is labeled $C(b)$ and not $C(a,b)$ because after intervention on $B$ there is no directed path from $A$ to $C$.
This labeling corresponds to the interpretation of missing edges in the graph in terms of the absence of individual level direct effects so that, for example, 
$C(a,b) = C(b)$ in Figure \ref{fig:swig-label}(d). \citet[\S 7]{thomas13swig} also discuss more general schemes that assume a time order, but also allow some missing edges to be interpreted at the individual level and others at the population (or distribution) level; in that paper ancestral labeling is termed `minimal labeling.'

Uniform labeling corresponds to the absence of any assumption regarding equality of potential outcomes (as random variables) across different interventions.\footnote{
Even in contexts where one assumes temporal and/or causal relationships between counterfactual random variables uniform labeling is useful when discussing more than one SWIG, since it makes clear which variable is in which SWIG; see the potential outcome calculus in \citet{malinsky19po} as an example.}
 In the potential outcome framework this would often appear somewhat unnatural. 
However, in this paper we will use this labeling to show that although we may wish to adopt the additional equalities between potential outcomes that are implied by the temporal and/or causal relationships, our results do not require these equalities. In addition, SWIGs with this labeling scheme are essentially isomorphic to the augmented decision diagrams proposed by \citet{dawid:2021}.  In particular, note that under the uniform labeling scheme the set of random variables appearing in two SWIGs $\G({ a})$ and $\G({ a^*})$, where $a,a^*\in \mathfrak{X}_A$  have no overlap; this will continue to hold when, below, we consider SWIGs $\G({ b})$ where we intervene on a (possibly empty) subset $B\subseteq A$.

\subsection{Distributional Consistency for SWIGs}

In order to relate passively observed distributions to those under intervention, we introduce a consistency assumption relating 
sets of counterfactual distributions. For this purpose we introduce the following notation:
\begin{align}
{\cal P}_A &\equiv \left\{\left.p(V({a})) \,\right|\, {a} \in \mathfrak{X}_A\right\},\\[6pt]
{\cal P}^{\tiny \subseteq}_A &\equiv \bigcup_{D\subseteq A} {\cal P}_D
\end{align}
Thus, ${\cal P}_A$ is the set of counterfactual distributions over $V$ that arise from all possible joint interventions setting the  variables in $A$ to a value ${a} \in \mathfrak{X}_A$. Likewise, ${\cal P}^{\subseteq}_A $ is the set of counterfactual distributions over $V$ resulting from all possible joint interventions on subsets $D$ of $A$; this includes the case $D=\emptyset$, corresponding to the observed distribution, so $p(V) \in {\cal P}^{\subseteq}_A $.

We make the following consistency assumption:\footnote{We thank Philip Dawid for pointing out an important omission in this definition in an earlier draft of this paper.}
\begin{definition}[Distributional Consistency for SWIGs]
The set of distributions ${\cal P}^{\subseteq}_A $ will be said to obey distributional consistency if, 
given  $B_i\in A$ and $C\subseteq  A\setminus \{B_i\}$, where $C$ may be empty, 
for all $y$, $b$, $c$:
\begin{align}\label{eq:consist-swig}
p(Y({ b},{ c})\!=\!y, B_i({ b},{ c}) \!=\! { b}) = p(Y({ c})\!=\!y, B_i({ c}) \!=\! { b}),
\end{align}
where $Y=V\setminus \{B_i\}$.
As a special case, if $C$ is empty then for all $y$, $b$:
\begin{align}\label{eq:consist-swig-simple}
p(Y({ b})\!=\!y , B_i({ b}) \!=\! { b}) = p(Y\!=\!y, B_i \!=\! { b}).
\end{align}

\end{definition}

The equalities (\ref{eq:consist-swig}) and (\ref{eq:consist-swig-simple}) simply state that the probability of the event $\{Y=y, B_i=b\}$, where here $B_i$ is the `natural' or (in Dawid's terminology) ITT variable, remains the same whether or not there is (subsequently) an intervention that  targets $B_i$ and sets it to $b$.

(\ref{eq:consist-swig-simple}) implies that $p(B_i(b)\!=\!b) = p(B_i\!=\!b)$,%
\footnote{In a standard potential outcome model, it would follow by definition of $B_i$, as the `natural value' of treatment, that $B_i(b)=B_i$.
Indeed, in the standard potential outcome approach there will not be a need to write $B_i(b)$. Readers familiar with the {\it do}-operator \citep{pearl95causal} should be aware that whereas in that theory intervention on a variable precludes observing the natural value, in the potential outcome theory, at least conceptually, we are supposing that the natural value could be observed and then, an instant later, we could intervene upon it without the natural value having any downstream causal effects.
\label{foot:natural-intervention}\par
Also note that the SWIG Local Markov Property (Definition \protect\ref{def:local-swig-mp}) will imply the stronger condition that $p(B_i(b^*)\!=\!b) = p(B_i\!=\!b)$ for all $b,b^* \in \mathfrak{X}_B$;
see Lemma \protect\ref{lem:margin-indep-future-fixed} below.
}
  and thus $p(Y(b)\!=\!y \mid B_i(b)\!=\!b) = p(Y\!=\!y \mid B_i\!=\!b)$. 
This has the interpretation that an intervention on $B_i$ setting it to $b$ is `ideal' in the sense that for the remaining variables $Y$, the intervention does not change the distribution of $Y$ given $B_i=b$. 
That $p(B_i(b)=b) = p(B_i=b)$ can be seen as following from the fact that $B_i$ and $B_i(b)$ represent, respectively, the natural value taken by $B_i$ in the absence of an intervention, and the natural value of $B_i$ immediately prior to an intervention. 

Under a standard potential outcome model 
that includes equalities between random variables 
(\ref{eq:consist-swig}) follows directly from the consistency assumption and recursive substitution:
\begin{align*}
B_i({ b},{ c}) \!=\! { b} \quad &\Rightarrow \quad B_i({ c}) \!=\! B_i({ b},{ c}) \!=\! { b}\\
\quad &\Rightarrow \quad Y({ b},{ c}) = Y(B_i({ c}),{ c}) = Y({ c}).
\end{align*} 
As with the discussion of labeling above, in a potential outcome theory it is natural to assume consistency at the level of random variables. Our motivation here for formulating consistency  via (\ref{eq:consist-swig}) as a relation between distributions is solely to make clear that we do not require the stronger assumption for our results.
However, proceeding in this way makes the notation more cumbersome since every potential outcome variable is labeled with every intervention.

Distributional consistency may also be formulated in terms of a dynamic regime. 
Let $g_i^*$ denote the dynamic regime\footnote{A dynamic regime is an intervention in which the value to which the variable is set is a function of the values taken by earlier variables.
With $g_i^*$ the earlier variable is the natural value that the variable would take on; see footnote \ref{foot:natural-intervention}.}
on $B$ which `intervenes' to set the intervention target to the `natural' value that the variable $B_i$ would take in the absence of an intervention. 
Let $V(g^*_i,c)$ be the set of potential outcomes that would arise under $g^*_i$ in conjunction with an intervention setting $C$ to $c$.
We may then re-express 
(\ref{eq:consist-swig}) as:
\begin{equation}\label{eq:swig-consist-dynamic}
p(V(g^*_i,c)) = p(V(c)).
\end{equation}
In words, in the context of an intervention setting $C$ to $c$, a dynamic regime which intervenes to set $B_i$ to the value that it would have taken anyway  has no effect on the distribution of $V$.
\footnote{Note that for (\ref{eq:swig-consist-dynamic}) to be equivalent to (\ref{eq:consist-swig}) we require that 
when $B_i(b,c) =b$, then $V(g^*_i,c) = V(b,c)$. This will hold if we assume: (i) that the values taken by variables that occur prior to the intervention on $B_i$ are unaffected by this intervention, and (ii) for variables that arise after the intervention, it makes no difference whether the value $b$ is imposed due to the dynamic regime $g^*_i$ and natural value 
$B_i(g^*,c) =b$, or due to a regime that uniformly imposes $b$. Both assumptions will hold if one views earlier variables, including the natural value $B_i(b,c)$, as having values that are determined prior to the decision to intervene on $B_i$. See \S \ref{sec:dawid-a} for more discussion of this point.}

Though the Distributional Consistency assumption involves a single variable $B_i$, repeated applications imply the same conclusion for a set $B$:

\begin{lemma}\label{lem:dist-consistency-vector}
If ${\cal P}^{\tiny \subseteq}_A$  obeys distributional consistency,   $B$, $C$ are disjoint subsets of $A$, where $C$ may be empty,
then 
for all $y$, $b$, $c$:
\begin{align}\label{eq:consist-swig-vector}
p(Y({ b},{ c})\!=\!y, B({ b},{ c}) \!=\! { b}) = p(Y({ c})\!=\!y, B({ c}) \!=\! { b}),
\end{align}
where $Y=V\setminus B$.
\end{lemma}

\noindent {\it Proof:} We prove this by induction on the size of $B$. The base case follows by definition of distributional consistency.
Let $B_i$ be a variable in $B$, and let $B_{-i} = B\setminus \{B_i\}$. 
\begin{align*}
\MoveEqLeft{p(Y({ b},{ c})=y, B(b,c)=b)}\\
 &= p(Y({ b_i, b_{-i}},{ c})=y, B_i(b_i,b_{-i},c)=b_i, B_{-i}(b_i,b_{-i},c)=b_{-i})\\
&= p(Y({b_{-i}},{ c})=y, B_i(b_{-i},c)=b_i, B_{-i}(b_{-i},c)=b_{-i})\\
&= p(Y({ c})=y, B_i(c)=b_i, B_{-i}(c)=b_{-i})\\
&= p(Y({ c})=y, B(c)=b).
\end{align*}
Here the second equality applies distributional consistency, taking `$C$' to be $B_{-i}\cup C$; the third applies the induction hypothesis,
taking `$Y$' to be $Y\cup \{B_i\}$ and `$B$' to be $B_{-i}$.\hfill$\Box$

The next Lemma relates equality of conditional distributions with and without an intervention on $B$.

\begin{lemma}\label{lem:cond-consistency} 
Suppose ${\cal P}^{\tiny \subseteq}_A$ obeys distributional consistency.
Let $B$, $C$ be disjoint subsets of $A$, where $C$ may be empty, and let $Y$, $W$ be disjoint subsets of $V\setminus B$. It then follows that:
\begin{align}\label{eq:consist-cond-swig}
p(Y({ b},{ c})\!=\!y \mid B({ b},{ c}) \!=\! { b}, W(b,c)\!=\!w) = p(Y({ c})\!=\!y \mid B({ c}) \!=\! { b}, W(c)\!=\!w).
\end{align}
\end{lemma}

\noindent{\it Proof:} This follows by applying Lemma \ref{lem:dist-consistency-vector} to $p(Y({ b},{ c}), B({ b},{ c}), W(b,c))$,
and $p(B({ b},{ c}), W(b,c))$. \hfill$\Box$\par
\medskip

In addition, we have the following:
\begin{lemma}\label{lem:remove-b-from-joint-not-a-fn}
Suppose ${\cal P}^{\tiny \subseteq}_A$ obeys distributional consistency, and let
$B$, $C$ be disjoint subsets of $A$, where $C$ may be empty. If $B\subseteq W \subseteq V$ and  $p(W({ b},{ c}))$ is not a function of 
${ b}$ then it follows from distributional consistency that $p(W({b},{c}))=p(W({c}))$.
\end{lemma}

\noindent{\it Proof:}
\begin{align*}
p(X_{W}(b,c)=w) &= p(X_{W\setminus B}(b,c)=w_{W\setminus B}, X_{B}(b,c)=w_{B})\\ 
&= p(X_{W\setminus B}(w_{B},c)=w_{W\setminus B}, X_{B}(w_{B},c)=w_{B})\\ 
&= p(X_{W\setminus B}(c)=w_{W\setminus B}, X_{B}(c)=w_{B})\\
&= p(X_{W}(c)=w).
\end{align*}
Here we use that $p(W({ b},{ c}))$ is not a function of 
${ b}$ in the second equality and distributional consistency via Lemma \ref{lem:dist-consistency-vector} in the third.
\hfill$\Box$
\\

Note that distributional consistency (\ref{eq:consist-swig}) does not imply the analogous result for conditional distributions. In particular, it is possible to have 
$B_i \in Y$,
$p(Y(b) \mid M(b))$ not be a function of $b$ and yet $p(Y(b) \mid M(b))\neq p(Y \mid M)$. This because even if $p(Y(b) \mid M(b))$ is not a function of $b$ both $p(Y(b),M(b))$ and $p(M(b))$ may still be functions of $b$, in which case there is no way to apply (\ref{eq:consist-swig}) to relate them to distributions in which $B$ is not intervened on.

However, when the conditioning set contains $B$ we have the following:

\begin{lemma}\label{lem:remove-b-from-conditional-not-a-fn}
Suppose ${\cal P}^{\tiny \subseteq}_A$ obeys distributional consistency, with $B$, $C$ disjoint subsets of $A$, where $C$ may be empty.  Further, let $Y$, $W$ be disjoint sets with $B\subseteq W$.
If $p(Y({ b},{ c})\mid W(b,c))$ is not a function of 
${ b}$ then
\begin{equation}\label{eq:remove-b-from-conditional}
p(Y({ b},{ c})\mid W(b,c)) = p(Y({ c})\mid W(c)).
\end{equation}
\end{lemma}

\noindent{\it Proof:}
\begin{align*}
\MoveEqLeft{p(X_{Y}(b,c)=y \mid X_{W}(b,c)=w)}\\
 &= p(X_{Y}(b,c)=y \mid X_{W\setminus B}(b,c)=w_{W\setminus B}, X_{B}(b,c)=w_{B} ) \\
  &= p(X_{Y}(w_{B},c)=y \mid X_{W\setminus B}(w_{B},c)=w_{W\setminus B}, X_{B}(w_{B},c)=w_{B} ) \\
   &= p(X_{Y}(c)=y \mid X_{W\setminus B}(c)=w_{W\setminus B}, X_{B}(c)=w_{B} )
\end{align*}
Here the second equality uses the fact that  $p(X_{Y}(b,c)=y \mid X_{W}(b,c)=w)$ is not a function of $b$, while the 
third follows from distributional consistency  
via Lemma \ref{lem:cond-consistency}.\hfill$\Box$

\subsection{Local Markov property defining the SWIG model}

Although we derive a SWIG graphically from the original DAG by node splitting, we will define the model by associating a local Markov property with the SWIG and the potential outcome distribution. The resulting model corresponds to the Finest Fully Randomized Causally Interpreted Structured Tree Graph (FFRCISTG) model of \citet{robins86new}; see  \citep[Appendix C]{thomas13swig}.
We will then derive the Markov property for the original DAG and the observed distribution from these by applying distributional consistency.

Given a DAG $\G$ with vertices $V=\{1,\ldots ,p\}$, we will use $\pa_{\G}(i)$ to indicate the (index) set of variables that are parents of $W_i$ in the original DAG $\G$, and let $\pre_{\prec}(i)$  indicate  $\{1,\ldots ,i-1\}$, the predecessors of $i$ under a total ordering $\prec$ that is consistent with the edges in $\G$. We will drop the subscript when the DAG or ordering is clear from context.

The SWIG local Markov property is defined on the set of distributions
${\cal P}_A \equiv \left\{\left. p(V({a}))\,\right| a \in \mathfrak{X}_{A}\right\}$, where $A\subseteq V$ is the maximal set of variables that may be intervened on; see \citet[\S1.2.4]{shpitser2021multivariate} and \citet{shpitser2021partially}.

\begin{definition}\label{def:local-swig-mp}
A set of potential outcome distributions ${\cal P}_A$ obeys {\em the SWIG ordered local Markov property 
for DAG $\G$ under $\prec$} if for all  $i \in V$, ${a} \in \frak{X}_A$, and  ${w} \in \frak{X}_{\pre_{\prec}(i)}$,
\begin{equation}\label{eq:swig:local-ind}
p\left(X_i({a}) \;\left|\; X_{\pre_{\prec}(i)}({a})\!=\! w \right.\right)
\end{equation}
is a function only of $a_{\pa_{\G}(i)\cap A}$ and $w_{\pa_{\G}(i)\setminus A}$.
\end{definition}

In words, (\ref{eq:swig:local-ind}) states that after intervening on $A$, the distribution of $X_i(a)$ given its predecessors depends solely on the values taken by intervention targets in $A$ that are parents of $i$, and by any other (random) variables that are parents of $i$ but that are not intervened on, and hence are not in $A$.\footnote{
In \citet[\S8, Def.44]{thomas13swig}, a weaker Markov property was stated that did not require that (\ref{eq:swig:local-ind}) is not a function of $a_{A\setminus \pa_{\G}(i)}$. 
This weaker condition does not imply the Markov property for the observed distribution (unless $A=V$). Consequently Propositions 45 and 46 and Theorem 65(c) in
\citet{thomas13swig} are incorrect. Correct reformulations are given below in Theorems \ref{thm:swig-obs-markov}, \ref{thm:swig-identify} and \ref{thm:swig-gb}.
}

Though the function of the local property is to define and characterize the potential outcome model, intuition may be gained by observing 
that the local property follows from d-separation applied to the SWIG $\G({ a})$.\footnote{This is as to be expected since d-separation encodes
the global property that is implied by the local property.} Specifically, the condition (\ref{eq:swig:local-ind}) corresponds to two sets of d-separations:
\medskip

\noindent{\bf d-separation from fixed nodes:}  That $p(X_i({a}) \,|\, X_{\pre(i)}({a}))$ does not depend on $a_{A \setminus \pa_{\G}(i)}$
is encoded in the SWIG $\G({ a})$ by the d-separation of  $X_i(a)$ from fixed nodes $a_j$ that correspond to vertices $A_j$ that are not parents of $X_i$ in $\G$ given the parents of $X_i(a)$ in $\G(a)$, both random and fixed; see  \citet{shpitser2021multivariate}, \citet{malinsky19po}, \citet{robins2018dseparation}.
  Specifically we have:
 \begin{align}\label{eq:d-sep-local-swig}
 X_i({a}) \indepd a_{A\setminus \pa(i)} \mid a_{A\cap \pa(i)}, X_{\pa(i)\setminus A}(a)
 \end{align}
 where here we used $\indepd$ to indicate d-separation\footnote{Note that in other papers \citep{shpitser2021multivariate,malinsky19po,thomas13swig} d-connection for SWIGs is defined such that fixed nodes may never occur as non-endpoint vertices on d-connecting paths. In those papers we never formally condition on fixed nodes.  Here, in (\ref{eq:d-sep-local-swig}) for the purpose of formulating the local property we formally include the fixed parents of $X_i(a)$ in the set that is (graphically) conditioned on. This is solely in order to make the development similar to the decision diagram approach we consider subsequently.%
 \label{foot:condition-on-fixed}}
 in the SWIG $\G(a)$ and use lower case letters, e.g.~$a_{A\setminus \pa(i)}$, to refer to fixed nodes.
  We may further decompose the set of fixed nodes $a_{A\setminus \pa(i)}$:
 \begin{align}\label{eq:d-sep-local-swig2}
 X_i({a}) \indepd \overbrace{\vphantom{X}a_{A\setminus \pre(i)}}^{\hbox{\tiny time order}},\;\;   
\overbrace{\vphantom{X}a_{(A\cap \pre(i))\setminus \pa(i)}}^{\hbox{\tiny causal Markov prop.}}\;\; \mid 
\underbrace{a_{A \cap \pa(i)},}_{\hbox{\tiny~~ fixed parents~~}}
\underbrace{X_{\pa(i) \setminus A}(a).}_{\hbox{\tiny ~~random parents~~}}
 \end{align}
The set of fixed nodes in $a_{A\setminus \pre(i)}$ correspond to interventions on variables that occur after $X_i$ and thus do not change 
$p(X_i({a}) \,|\, X_{\pre(i)}({a}))$. Likewise, the effects of the fixed nodes in $a_{(A\cap \pre(i))\setminus \pa(i)}$ are screened off by the random and fixed nodes that are parents of $X_i(a)$.
\medskip

\noindent{\bf d-separation from random nodes:}  That $p(X_i({a}) \mid X_{\pre(i)}({a})\!=\! w )$ does not depend on $w_{\pre(i) \setminus (\pa_{\G}(i)\setminus A}$ is encoded in $\G(a)$ by the d-separation of $X_i({a})$ from $X_{\pre(i) \setminus (\pa(i)\setminus A)}(a)$ conditioning on 
the parents of $X_i(a)$ in $\G(a)$, both random and fixed:
 \begin{align}\label{eq:d-sep-local-swig3}
 X_i({a}) \indepd X_{\pre(i) \setminus (\pa(i)\setminus A)}(a) \;\mid\; a_{A\cap \pa(i)}, X_{\pa(i)\setminus A}(a).
 \end{align}
The random vertices $X_{\pre(i) \setminus (\pa(i)\setminus A)}(a)$ may be further decomposed:
 \begin{align}\label{eq:d-sep-local-swig4}
 X_i({a}) \indepd \;\;\overbrace{X_{\pre(i)\setminus \pa(i)}(a),}^{\hbox{\tiny assoc.~Markov prop.}} \quad \overbrace{X_{\pa(i) \cap A}(a)}^{\hbox{\tiny ignorability}} \;\;\mid 
 \underbrace{a_{A \cap \pa(i)},}_{\hbox{\tiny~~ fixed parents~~}}
\underbrace{X_{\pa(i) \setminus A}(a).}_{\hbox{\tiny ~~random parents~~}}
 \end{align}
The d-separation of $X_i(a)$ from nodes representing the natural value of variables that are in $A$ and parents of $X_i$ in $\G$ corresponds to ignorability. On the other hand, the d-separation of $X_i(a)$ from variables that are predecessors, but not parents, of $X_i$ in $\G$ can be regarded as an associational Markov property.

\subsection{Example}

The d-separations given by (\ref{eq:d-sep-local-swig2}) and (\ref{eq:d-sep-local-swig4}) can obviously be stated as a single graphical condition for each random vertex $V_i(a)$ in $\G(a)$. In Tables \ref{tab:swig-factor} and \ref{tab:swig-dsep} we give the SWIG local Markov property corresponding to the SWIG $\G({\bf x}) \equiv \G(x_1,x_2)$ shown in Figure \ref{fig:pearl-wrong-example2}(b),\footnote{\citet[Figure 15]{dawid:2021} gives the corresponding SWIG under ancestral labeling. In his discussion of this example \citet[Figure 15]{dawid:2021}, two conditions are stated as supporting g-computation. The first of these is correct, but the second should be $Y(x_0,x_1) \indep X_0$, not $Z(x_0) \indep X_0$.} under the ordering $(H,X_0,Z,X_1,Y)$:
 Table \ref{tab:swig-factor} in terms of factorization; Table \ref{tab:swig-dsep} via d-separation.
Note that for each $V_i$, the number  of arguments on which $p(V_i({\bf x}) \,|\, V_{\pre(i)}({\bf x}))$ depends corresponds exactly to the number of parents (random and fixed) of the corresponding random variable in $\G({\bf x})$ in Figure \ref{fig:pearl-wrong-example2}(b); zero for $H(\bf{x})$, $X_0(\bf{x})$; two for $Z(\bf{x})$, $X_1(\bf{x})$ and $Y(\bf{x})$. This is also the number of terms listed to the right of the conditioning bar in Table \ref{tab:swig-dsep}.
Here, as elsewhere in this paper, we use the uniform labeling because we wish to emphasize that our results do not require any equalities between random variables. 

{\begin{table}
\centering
 \begin{tabular}{r@{\extracolsep{1pt}}l}
\multicolumn{2}{c}{local Markov property for $\G(x_1,x_2)$ via factorization terms}\\[5pt]
\hline\\[-5pt]
 $p(H({\textcolor{red}{x_0}},{\textcolor{red}{x_1}}))$ &\\[2pt]
 $p(X_0({\textcolor{red}{x_0}},{\textcolor{red}{x_1}})\;$& $\mid {\textcolor{red}{H(x_0,x_1)}})$\\[2pt]
 $p(Z({x_0},{\textcolor{red}{x_1}})\;$& $\mid {H({x_0},\textcolor{red}{x_1})},{\textcolor{red}{X_0(x_0,x_1)}})\;$\\[2pt]
 $p(X_1(\textcolor{red}{x_0},{\textcolor{red}{x_1}})\;$& $\mid {H(\textcolor{red}{x_0},\textcolor{red}{x_1})},{\textcolor{red}{X_0(x_0,x_1)}}, {Z(\textcolor{red}{x_0},\textcolor{red}{x_1})})$\\[2pt]
$p(Y(\textcolor{red}{x_0},x_1)\;$& $\mid {\textcolor{red}{H(x_0,x_1)}}, {\textcolor{red}{X_0(x_0,x_1)}}, {Z(\textcolor{red}{x_0},{x_1})}, {\textcolor{red}{X_1(x_0,x_1)}})$\\[2pt]
\end{tabular}
\caption{Defining properties for the SWIG $\G({\bf x})$ in Figure \ref{fig:pearl-wrong-example2}(b), expressed via factorization. Arguments in $p(V_i({\bf x}) \,|\, V_{\pre(i)}({\bf x}))$ on which this term does not depend are colored red. Note that the arguments in $V_{\pre(i)}({\bf x}))$ on which the term depends, correspond to the parents of $V_i({\bf x})$ in $\G({\bf x})$; these are written in black. For example, for the term corresponding to $V_i=Y$, the arguments are $x_1$ and $Z({\bf x})$, and these are the parents of $Y({\bf x})$ in $\G({\bf x})$.
\label{tab:swig-factor}}
\end{table}

\begin{table}
\centering
\begin{tabular}{r@{\extracolsep{2pt}}c@{\extracolsep{2pt}}l}
\multicolumn{3}{c}{local Markov property for $\G(x_1,x_2)$ via d-separation}\\[5pt]
\hline\\[-5pt]
 $H(x_0,x_1) \indepd$ &$ x_0,x_1$\\[2pt]
 $X_0(x_0,x_1) \indepd$ &$ H(x_0,x_1), x_0,x_1$ \\[2pt]
 $Z(x_0,x_1) \indepd$ &$ X_0(x_0,x_1), x_1 $ &$\mid H(x_0,x_1), x_0 $\\[2pt]
  $X_1(x_0,x_1) \indepd$ &$ X_0(x_0,x_1), x_0, x_1 $ &$\mid H(x_0,x_1), Z(x_0,x_1) $\\[2pt]
  $Y(x_0,x_1) \indepd$ &$ X_0(x_0,x_1), X_1(x_0,x_1), H(x_0,x_1), x_0 $ &$\mid Z(x_0,x_1), x_1 $ \\[20pt]
 \end{tabular}
\caption{The d-separation relations corresponding to the SWIG local Markov property in the SWIG $\G(x_1,x_2)$ in Figure \ref{fig:pearl-wrong-example2}(b).
Here $x_1$ and $x_2$ refer to the fixed nodes and $\indepd$ indicates d-separation in the SWIG; see also footnote \protect\ref{foot:condition-on-fixed} regarding
the formal inclusion of fixed nodes on the RHS of the conditioning bar.
\label{tab:swig-dsep}}
\end{table}

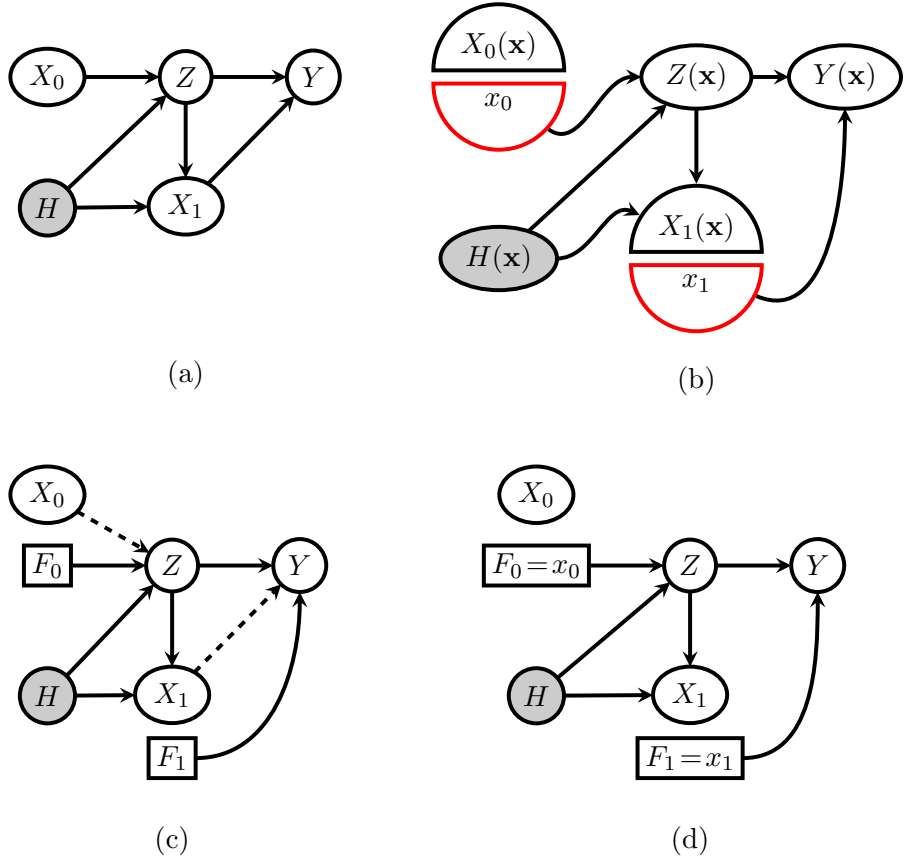
\begin{figure}
\begin{center}
\begin{tikzpicture}[
pre/.style={->,>=stealth,semithick,blue,line width = 1pt}]
\tikzset{line width=1.25pt, outer sep=0pt,
         ell/.style={draw, inner sep=3pt,
          line width=1.5pt}, >=stealth,
         swig vsplit={gap=5pt,
         inner line width right=0pt, line color right=red},
         ell2/.style={draw,fill opacity=0, text opacity=1}};
\begin{scope}
\node[ell, name=a,shape=ellipse,style={draw}] {
$X_0$
};
\node[ell, name=z,shape=circle,style={draw},right=1cm of a]
{$Z$
};
\node[ell, name=h,fill=black!20,shape=circle,style={draw},below =1cm of a]
{$H$
};
\node[ell, name=b,shape=ellipse,style={draw},below =1cm of z]
{$X_1$
};
\node[ell, name=y,shape=circle,style={draw},right =1cm of z]
{$Y$
};
\node[below =3.3cm of z]{(a)};
\draw[ell,->] (a) to (z);
\draw[ell,->] (z) to (y);
\draw[ell,->] (b) to (y);
\draw[ell,->] (z) to (b);
\draw[ell,->] (h) to (b);
\draw[ell,->] (h) to (z);
\end{scope}
\begin{scope}[xshift=6cm]
\node[ell,name=a,
shape=swig hsplit, swig hsplit={line color lower=red,gap=5pt}, ell2]{
        \nodepart{upper}{$X_0({\bf x})$}
        \nodepart{lower}{$x_0$}   };
\node[ell, name=z,shape=ellipse,style={draw},right=1cm of a]
{$Z({\bf x})$
};
\node[ell, name=h,fill=black!20,shape=ellipse,style={draw},below =1cm of a]
{$H({\bf x})$
};
\node[ell,name=b, below=1cm of z,
shape=swig hsplit, swig hsplit={line color lower=red,gap=5pt}, ell2]{
        \nodepart{upper}{$X_1({\bf x})$}
        \nodepart{lower}{$x_1$}   };
\node[ell, name=y,shape=ellipse,style={draw},right =0.5cm of z]
{$Y({\bf x})$
};
\node[below =3.3cm of z]{(b)};
\draw[ell,->] (a) to[out=-40,in=180] (z);
\draw[ell,->] (z) to (y);
\draw[ell,->] (b) to[out=-25,in=270] (y);
\draw[ell,->] (z) to (b);
\draw[ell,->] (h) to[out=0,in=150] (b);
\draw[ell,->] (h) to (z);
\end{scope}
\begin{scope}[yshift=-6.5cm]
\node[ell,name=fa,shape=rectangle,style={draw}] {
$F_0$
};
\node[ell,name=a,above=3mm of fa, shape=ellipse,style={draw}] {
$X_0$
};
\node[ell, name=z,shape=circle,style={draw},right=1cm of fa]
{$Z$
};
\node[ell, name=h,fill=black!20,shape=circle,style={draw},below =1.1cm of fa]
{$H$
};
\node[ell, name=b,shape=ellipse,style={draw},below =1cm of z]
{$X_1$
};
\node[ell,name=fb,below=2mm of b,shape=rectangle,style={draw}] {
$F_1$
};
\node[ell, name=y,shape=circle,style={draw},right =1cm of z]
{$Y$
};
\draw[ell,dashed,->] (a) to (z);
\draw[ell,->] (fa) to (z);
\draw[ell,->] (z) to (y);
\draw[ell,dashed,->] (b) to (y);
\draw[ell,->] (fb) to[out=0,in=270] (y);
\draw[ell,->] (z) to (b);
\draw[ell,->] (h) to (b);
\draw[ell,->] (h) to (z);
\node[below =3cm of z]{(c)};
\end{scope}
\begin{scope}[xshift=6.5cm,yshift=-6.5cm]
\node[ell,name=fa,shape=rectangle,style={draw}] {
$F_0\!=\!x_0$
};
\node[ell,name=a,above=3mm of fa, shape=ellipse,style={draw}] {
$X_0$
};
\node[ell, name=z,shape=circle,style={draw},right=1cm of fa]
{$Z$
};
\node[ell, name=h,fill=black!20,shape=circle,style={draw},below =1.1cm of fa]
{$H$
};
\node[ell, name=b,shape=ellipse,style={draw},below =1cm of z]
{$X_1$
};
\node[ell,name=fb,below=2mm of b,shape=rectangle,style={draw}] {
$F_1\!=\!x_1$
};
\node[ell, name=y,shape=circle,style={draw},right =1cm of z]
{$Y$
};
\node[below =3cm of z]{(d)};
\draw[ell,->] (fa) to (z);
\draw[ell,->] (z) to (y);
\draw[ell,->] (fb) to[out=0,in=270] (y);
\draw[ell,->] (z) to (b);
\draw[ell,->] (h) to (b);
\draw[ell,->] (h) to (z);
\end{scope}
\end{tikzpicture}
\end{center}
\caption{(a) The DAG $\G$ originally considered in Ex.~11.3.3, Fig.~11.12 in \protect{\cite[p.353]{pearl:2009}}, here $H$ is unobserved; (b)
the SWIG ${\G}({\bf x})$ under uniform labeling, ${\bf x}\equiv (x_1,x_2)$; Figure 15 in \cite{dawid:2021} shows the SWIG with ancestral labeling; (c) the reformulated augmented graph $\G^*$ (this corresponds to Dawid's ITT DAG $\G^*$ shown in 
Figure 13 of \protect{\cite[p.62]{dawid:2021}} after marginalizing $X_0$, $X_1$ and then removing $^*$ from the ITT variables);
(d) the resulting graph under the regime $F_0=x_0$, $F_1=x_1$; this is the graph that encodes the reformulated Markov property.
\label{fig:pearl-wrong-example2}}
\end{figure}

\subsection{Consequences of the Local Markov property}

Under distributional consistency it follows from the SWIG local Markov property that whether or not interventions in the future occur has no effect on the distribution of prior variables.

\begin{lemma}\label{lem:margin-indep-future-fixed}
 If ${\cal P}^{\subseteq}_A$ obeys distributional consistency and ${\cal P}_A$  obeys the 
SWIG ordered local Markov property for DAG $\G$ under $\prec$ then for all $k \in V$ and $a\in {\mathfrak{X}}_A$:
\begin{equation}\label{eq:no-future-effect}
p(X_1(a),\ldots ,X_k(a)) = p(X_1(a_{\pre(k)\cap A}),\ldots ,X_k(a_{\pre(k)\cap A})).
\end{equation}
\end{lemma}

\noindent{\it Proof:}  First observe that 
since
\begin{align*}
p(X_1(a),\ldots ,X_k(a)) = \prod_{i=1}^k p(X_i({ a}) \mid X_{\pre(i)}({ a}) ),
\end{align*}
and the local Markov property implies that $p(X_i({ a}) \mid X_{\pre(i)}({ a}) )$ does not depend on $a_{A\setminus \pre(i)}$, 
it follows that $p(X_1(a),\ldots ,X_k(a))$ does not depend on $a_{A\setminus \pre(k)}$.

We now prove the claim by reverse induction on the ordering of the vertices in $V$.\par
For the base case suppose $k$ is the maximal vertex in $V$. If $k\notin A$ then (\ref{eq:no-future-effect}) holds trivially
since $A = \pre(k)\cap A$.
If $k\in A$ then since $k\notin \pre(k)$, $p(X_1(a),\ldots ,X_k(a))$ does not depend on $a_k$ and thus, by 
Lemma \ref{lem:remove-b-from-joint-not-a-fn}, $p(X_1(a),\ldots ,X_k(a))=p(X_1(a_{A\setminus\{k\}}),\ldots ,X_k(a_{A\setminus\{k\}}))$.\par
Our inductive hypothesis is that (\ref{eq:no-future-effect}) holds for $k=j+1$, so that  
\begin{equation*}
p(X_1(a),\ldots ,X_{j+1}(a)) = p(X_1(a_{\pre(j+1)\cap A}),\ldots ,X_{j+1}(a_{\pre(j+1)\cap A})).
\end{equation*}
Summing both sides over $x_{j+1}$ we obtain:
\begin{equation}\label{eq:inductive-case-j}
p(X_1(a),\ldots ,X_{j}(a)) = p(X_1(a_{\pre(j+1)\cap A}),\ldots ,X_{j}(a_{\pre(j+1)\cap A})).
\end{equation}
If $j \notin A$ then (\ref{eq:inductive-case-j}) establishes the claim since $\pre(j+1)\cap A = \pre(j)\cap A$. If $j \in A$ then note that we have already established above that the LHS of (\ref{eq:inductive-case-j})
is not a function of $a_j$. Consequently, the RHS is also not a function of $a_j$. It then follows from Lemma \ref{lem:remove-b-from-joint-not-a-fn} that
\begin{align*}
\MoveEqLeft{p(X_1(a_{\pre(j+1)\cap A}),\ldots ,X_{j}(a_{\pre(j+1)\cap A}))}\\
&= p(X_1(a_{\pre(j)\cap A}),\ldots ,X_{j}(a_{\pre(j)\cap A})).
\end{align*}
This completes the proof. \hfill$\Box$\par
\bigskip

The next Lemma gives a simple characterization of the consequences of the SWIG local Markov property in conjunction with distributional consistency.

\begin{lemma}\label{lem:kernel-of-swig-lmp}
 If ${\cal P}^{\subseteq}_A$ obeys distributional consistency and ${\cal P}_A$  obeys the 
SWIG ordered local Markov property for DAG $\G$ under $\prec$ then:
\begin{align}
\MoveEqLeft{p(X_i({ a}) \mid X_{\pre(i)}({ a}) )}\label{eq:swig-lmp-lhs}\\
&= p(X_i({ a}_{\pre(i)\cap A}) \mid X_{\pre(i)}({ a}_{\pre(i)\cap A}) ) \label{eq:drop-successors-from-a}  \\
&= p(X_i({ a}_{\pa(i)\cap A}) \mid X_{\pre(i)}({ a}_{\pa(i)\cap A}) ) \label{eq:drop-predecessors-from-a} \\
&= p(X_i({ a}_{\pa(i)\cap A}) \mid X_{\pa(i)}({ a}_{\pa(i)\cap A}) ) \label{eq:drop-predecessors-from-x} \\
&= p(X_i({ a}_{\pa(i)\cap A}) \mid X_{\pa(i)\setminus A}({ a}_{\pa(i)\cap A}) ). \label{eq:drop-parents-from-x} 
\end{align}
\end{lemma}

Since  the SWIG local Markov property (\ref{eq:swig:local-ind}) states that (\ref{eq:swig-lmp-lhs}) is not a function of
$a_{A \setminus \pa_{\G}(i)}$, the
equality of (\ref{eq:swig-lmp-lhs}) and (\ref{eq:drop-predecessors-from-a}) 
may appear to follow immediately.
However, as noted in the discussion prior to Lemma \ref{lem:remove-b-from-conditional-not-a-fn}, the fact that a counterfactual conditional distribution $p(Y(a_j)\,|\,W(a_j))$ does not depend on the specific value, $a_j$, of an intervention on $A_j$ does {\em not} imply that
$p(Y(a_j)\,|\,W(a_j)) = p(Y\,|\,W)$.
\medskip

\noindent{\it Proof:} Here (\ref{eq:drop-successors-from-a}) follows since by Lemma \ref{lem:margin-indep-future-fixed}
\begin{align*}
p(X_i({ a}), X_{\pre(i)}({ a}) ) 
&= p(X_i({ a}_{\pre(i)\cap A}), X_{\pre(i)}({ a}_{\pre(i)\cap A}) ).
\end{align*}
(\ref{eq:drop-predecessors-from-a}) follows from Definition \ref{def:local-swig-mp} and Lemma \ref{lem:remove-b-from-conditional-not-a-fn}.  Finally, (\ref{eq:drop-predecessors-from-x}) and (\ref{eq:drop-parents-from-x}) follow from the SWIG local Markov property via (\ref{eq:remove-b-from-conditional}) 
since $p(X_i({ a}) \mid X_{\pre(i)}({ a})\!=\! x_{\pre(i)} )$ does not depend on $x_{\pre(i)\setminus (\pa(i) \setminus A)} = 
(x_{\pre(i)\setminus \pa(i)}, x_{\pa(i)\cap A})$.  \hfill$\Box$\par

\subsection{Markov property for the observed distribution}

We now show that distributional consistency together with the SWIG local Markov property implies the usual local Markov property 
\citep{lau:indmark} for the observed distribution. 

\begin{theorem}\label{thm:swig-obs-markov}
 If ${\cal P}^{\subseteq}_A$ obeys distributional consistency and ${\cal P}_A$  obeys the 
SWIG ordered local Markov property for $\G$  and $\prec$ then $p(V)$ obeys the usual ordered local Markov property w.r.t. $\G$ and $\prec$.
\end{theorem}

\noindent{\it Proof:}  Let $v^* \in \mathfrak{X}_{\pre(i)}$.
\begin{align}
\MoveEqLeft{p(X_i=v \mid X_{\pre(i)}=v^*)}\nonumber\\
&= p(X_i(v^*_{\pre(i)\cap A })=v \mid X_{\pre(i)}(v^*_{\pre(i)\cap A})\!=\!v^*)\nonumber\\
&= p(X_i(v^*_{\pa(i) \cap A})=v \mid X_{\pa(i)\setminus A}(v^*_{\pa(i) \cap A})\!=\!v^*_{\pa(i)\setminus A}).
\end{align}
Here the first equality follows from distributional consistency via Lemma \ref{lem:cond-consistency}.
The second follows directly from the equality of (\ref{eq:drop-successors-from-a}) and (\ref{eq:drop-parents-from-x}) in Lemma \ref{lem:kernel-of-swig-lmp}.
Since the last line is not a function of $v^*_{\pre(i)\setminus \pa(i)}$, the ordered local Markov property for the DAG holds.
\hfill$\Box$
\medskip

\subsubsection{Discussion of relation to Dawid}

Dawid takes the reverse approach to ours: he proposes additional extended Markovian conditions that, when added to the usual Markov property for the observable law will imply the Markov property for his extended graph. However, as we describe in detail below, our approach appears to be  simpler in that, given distributional consistency, it requires only one property per variable, giving $|V|$ constraints in total; in contrast, Dawid requires one property for every observed variable in $V$, together with two additional properties for each intervention target in $A$ for a total of ($|V|+2|A|$).

In addition, our approach captures context specific independences, corresponding to `dashed' edges in Dawid's diagrams; further these are not captured directly in Dawid's A+B formulation.
 We show that by restating the SWIG local property in Dawid's notation, we are able to provide a characterization of the (extended) Markov properties for the augmented graph and the original graph that also requires only one constraint per variable, plus distributional consistency.

It is the case that Dawid incorporates distributional consistency into his defining independences, whereas we state it as a separate property that precedes the definition of the model. However, as we have shown above, distributional consistency may be seen as a tautologous property the truth of which is implicit in the notion of an ideal intervention: distributional consistency states that if $B$ would naturally take the value $b$, then an ideal intervention that would set $B$ to $b$ has no effect on the distribution of (all) the other variables. For this reason, we believe it is natural to distinguish consistency from the other properties being used to define the model.

However, in the spirit of Dawid's approach, in Appendix \ref{app:dawid-swig} we show that 
if ${\cal P}^{\subseteq}_A$ obeys distributional consistency, then 
${\cal P}_A$ will obey the SWIG local Markov property corresponding to ${\G}$ if:  
(i) $p(V)$ is positive and obeys the (ordinary) local Markov property for the graph $\G$;
(ii) ${\cal P}_A$ obeys the SWIG local Markov property corresponding to $\overline{\G}$, a complete supergraph of $\G$. This formulation requires $2|V|$ restrictions.

\subsection{Identification of the potential outcome distribution $p(V(a))$ from $p(V)$}

We show that it follows from the SWIG local Markov property that $p(V(a))$ is identified given the distribution over the observables provided that the relevant conditional distributions are identified from the distribution of the observables.

\begin{theorem}\label{thm:swig-identify}
Suppose that ${\cal P}^{\subseteq}_A$ obeys distributional consistency and ${\cal P}_A$  obeys the 
SWIG ordered local Markov property for $\G$  and $\prec$.
Let $a\in  \mathfrak{X}_{A}$ be an assignment to the intervention targets in $A$, and let $v \in \mathfrak{X}_{V}$. Then for all $i$:
\begin{align}
\MoveEqLeft{p(X_i({ a})=v_i \mid X_{\pre(i)}({ a})= v_{\pre(i)})}\nonumber\\
&= p(X_i=v_i \mid  X_{\pa(i) \setminus A} = v_{\pa(i) \setminus A}, X_{\pa(i)\cap A}=a_{\pa(i)\cap A}).\label{eq:modularity-swig}
\end{align}
Consequently, $p(V({a}))$ is identified from $p(V)$ and obeys d-separation in the SWIG $\G({a})$,
whenever the conditional distributions on the RHS of {\rm (\ref{eq:modularity-swig})} are
identified by $p(V)$.
\end{theorem}

\noindent The equality (\ref{eq:modularity-swig}) here corresponds to the property referred to as `modularity' in \citep{thomas13swig}; this is also an instance
of the extended g-formula of \citet{robins86new,robins:effects:2004}.\\

\noindent{\it Proof:} 
\begin{align}
\MoveEqLeft{p(X_i(a)\!=\! v_i \mid X_{\pre(i)}(a)\!=\! v_{\pre(i)})
}\nonumber\\
&=  p(X_i({ a}_{\pa(i)\cap A})\!=\! v_i \mid X_{\pa(i)\cap A}({ a}_{\pa(i)\cap A})\!=\! v_{\pa(i)\cap A},\nonumber\\
&\kern130pt X_{\pa(i)\setminus A}({ a}_{\pa(i)\cap A})\!=\! v_{\pa(i)\setminus A} ) \nonumber\\[4pt]
&=  p(X_i({ a}_{\pa(i)\cap A})\!=\! v_i \mid X_{\pa(i)\cap A}({ a}_{\pa(i)\cap A})\!=\! a_{\pa(i)\cap A},\nonumber\\
&\kern130pt X_{\pa(i)\setminus A}({ a}_{\pa(i)\cap A})\!=\! v_{\pa(i)\setminus A} ) \nonumber\\[4pt]
&=  p(X_i\!=\! v_i \mid X_{\pa(i)\cap A}\!=\! a_{\pa(i)\cap A}, X_{\pa(i)\setminus A}\!=\! v_{\pa(i)\setminus A} ).
\end{align}
Here the first equality follows from the equality of (\ref{eq:swig-lmp-lhs}) and (\ref{eq:drop-predecessors-from-x});  the second follows from
the equality of (\ref{eq:drop-predecessors-from-x}) and (\ref{eq:drop-parents-from-x}); the third from distributional consistency via (\ref{eq:consist-cond-swig}).
\hfill$\Box$

\subsection{Distributions resulting from fewer interventions}

Finally we show that if $p(V(a))$ obeys the SWIG local Markov property for $\G$ and distributional consistency, then
if we intervene on $B \subset A$, the resulting distribution $p(V(b))$ will obey the SWIG local Markov property for $\G$ with respect to this reduced set of intervention targets. The two previous Theorems can be seen as the special case in which $B=\emptyset$.

\begin{theorem}\label{thm:swig-gb}
Suppose that ${\cal P}^{\subseteq}_A$ obeys distributional consistency and ${\cal P}_A$  obeys the 
SWIG ordered local Markov property for $\G$  and $\prec$.
Let ${b}$ be an assignment to the intervention targets in $B\subseteq A$, and let $v\in \mathfrak{X}_V$. Then for all $i$:
\begin{align}
\MoveEqLeft{p(X_i({ b}) =v_i \mid X_{\pre(i)}({b}) = v_{\pre(i)}) 
}
\nonumber\\
&= p(X_i=v_i \mid X_{\pa(i)\cap B}=b_{\pa(i)\cap B}, X_{\pa(i) \setminus B} = v_{\pa(i) \setminus B}).\label{eq:modularity-for-b-swig}
\end{align}
Consequently, every $p(V({b})) \in {\cal P}_B$ 
obeys the Markov property for the SWIG $\G({b})$ and is identified whenever the conditional distributions on the RHS of 
{\rm (\ref{eq:modularity-for-b-swig})}
are identified by $p(V)$.
\end{theorem}

\noindent{\it Proof:} 
\begin{align}
\MoveEqLeft{p(X_i({ b}) =v_i \mid X_{\pre(i)}({ b})\!=\! v_{\pre(i)}) }\nonumber\\
&= p(X_i({ b_{\pre(i)\cap B}}) =v_i \mid X_{\pre(i)}({b_{\pre(i)\cap B}})\!=\! v_{\pre(i)}) \nonumber\\[2pt]
&= p(X_i({ b_{\pre(i)\cap B}}, v_{\pre(i)\cap (A\setminus B)}) =v_i \mid\nonumber \\[-3pt]
&\kern120pt X_{\pre(i)}({b_{\pre(i)\cap B}}, v_{\pre(i)\cap (A\setminus B)})\!=\! v_{\pre(i)}) \nonumber\\[2pt]
&= p(X_i =v_i \mid X_{\pa(i)\setminus A}\!=\! v_{\pa(i)\setminus A},
X_{\pa(i)\cap B}\!=\! b_{\pa(i)\cap B},\nonumber\\[-3pt]
&\kern170pt X_{\pa(i)\cap ( A\setminus B) }\!=\! v_{\pa(i)\cap (A\setminus B)}) \nonumber\\[2pt]
&= p(X_i =v_i \mid X_{\pa(i)\setminus B}\!=\! v_{\pa(i)\setminus B},
X_{\pa(i)\cap B}\!=\! b_{\pa(i)\cap B}) \nonumber
\end{align}

Here the first equality is by Lemma \ref{lem:remove-b-from-joint-not-a-fn}; the second is distributional consistency via Lemma \ref{lem:cond-consistency}; 
the third follows from Theorem \ref{thm:swig-identify} applied to $\G(a)$; the fourth is a simplification.\hfill$\Box$
\begin{figure}
\centering
\begin{tikzpicture}[scale=0.9]
\tikzset{line width=1.25pt, outer sep=0pt,
         ell/.style={draw, inner sep=2pt,
          line width=1.5pt}, >=stealth,
         swig vsplit={gap=5pt,
         inner line width right=0pt, line color right=red},
         ell2/.style={draw,fill opacity=0, text opacity=1}}; 
\begin{scope}
\node[name=t,ell,draw,shape=ellipse]{$T^*$};
\node[name=y,ell,draw,shape=ellipse,right=15mm of t]{$Y$};
\draw[ell,->](t) to (y);
\node[left=3mm of t]{(a)};
\end{scope}
\begin{scope}[xshift=6cm,yshift=5cm]
\node[name=t,ell,shape=ellipse, draw=none,opacity=0]{\phantom{$T^*$}};
\node[name=tplus,ell,draw,shape=ellipse, right = 10mm of t]{$T$};
\node[name=f,ell,draw,shape=rectangle, inner sep=4pt, above = 5mm of tplus]{$F_T$};
\node[name=y,ell,draw,shape=ellipse,right=8mm of tplus]{$Y$};
\draw[ell,->](tplus) to (y);
\draw[ell,->](f) to (tplus);
\node[left=3mm of t]{(c)};
\end{scope}
\begin{scope}[yshift=-2.5cm,xshift=5mm]         
\node[name=t, shape=swig vsplit, ell2]{
        \nodepart{left}{$T^*$}
        \nodepart{right}{$t$}};
\node[name=y,ell,draw,shape=ellipse,right=5mm of t]{$Y(t)$};
\draw[ell,->](t) to (y);
\node[left=3mm of t]{(b)};
\end{scope}
\begin{scope}[xshift=6cm, yshift=2.5cm]
\node[name=t,ell,draw,shape=ellipse]{$T^*$};
\node[name=tplus,ell,draw,shape=ellipse, right = 10mm of t]{$T$};
\node[name=f,ell,draw,shape=rectangle, inner sep=4pt, above = 5mm of tplus]{$F_T$};
\node[name=y,ell,draw,shape=ellipse,right=8mm of tplus]{$Y$};
\draw[ell,->, red, dashed, line width=1.75pt](t) to (tplus);
\draw[ell,->](tplus) to (y);
\draw[ell, red, ->](f) to (tplus);
\node[left=3mm of t]{(d)};
\end{scope}
\begin{scope}[xshift=6cm]
\node[name=t,ell,draw,shape=ellipse]{$T^*$};
\node[name=tplus,ell,draw,shape=ellipse, right = 10mm of t]{$T$};
\node[name=f,ell,draw,shape=rectangle, inner sep=4pt, above = 5mm of tplus]{$F_T=\emptyset$};
\node[name=y,ell,draw,shape=ellipse,right=8mm of tplus]{$Y$};
\draw[ell,->,red, line width=1.75pt](t) to (tplus);
\draw[ell,->](tplus) to (y);
\node[left=3mm of t]{(e)};
\end{scope}
\begin{scope}[xshift=6cm,yshift=-2.5cm]
\node[name=t,ell,draw,shape=ellipse]{$T^*$};
\node[name=tplus,ell,draw,shape=ellipse, right = 5mm of t]{$T=t$};
\node[name=f,ell,draw,shape=rectangle, inner sep=4pt, above = 5mm of tplus]{$F_T\!=\!t$};
\node[name=y,ell,draw,shape=ellipse,right=5mm of tplus]{$Y$};
\draw[ell,->,red, line width=1.75pt](f) to (tplus);
\draw[ell,->](tplus) to (y);
\node[left=3mm of t]{(f)};
\end{scope}
\begin{scope}[xshift=6cm,yshift=-5cm]
\node[name=t,ell,draw,shape=ellipse]{$T^*$};
\node[name=tplus,ell,draw=none,opacity=0,shape=ellipse, right = 10mm of t]{\phantom{$T$}};
\node[name=f,ell,draw,shape=rectangle, inner sep=4pt, above = 5mm of tplus]{$F_T$};
\node[name=y,ell,draw,shape=ellipse,right=8mm of tplus]{$Y$};
\draw[ell,->,dashed, line width=1.75pt](t) to (y);
\draw[ell,->](f) to (y);
\node[left=3mm of t]{(g)};
\end{scope}
\end{tikzpicture}
\caption{ The simplest case of a single treatment $T$ and outcome $Y$ in the absence of confounding.
{\small (a) DAG $\G$ representing the observed joint distribution $p(T^*,Y)$; (b) SWIG ${\G}(t)$ corresponding to ${\G}$ representing $p(T^*,Y(t))$; (c) Dawid's augmented DAG representing the set of kernels $p(Y,T \mid F_T)$ where $F_T$ is a regime indicator;
(d) Dawid's augmented DAG with intention-to-treat (ITT) variables, representing the kernels $p(Y,T^*,T \mid F_T)$ where $F_T$ is a regime indicator; the dashed edge indicates that the edge between $T^*$ and $T$ is absent in the interventional regime, while the red edges indicate deterministic relationships;
 (e) the ITT augmented graph representing the observational regime $p(T^*,T, Y\mid F_T=\emptyset) = p(T^*,T, Y)$ under which $T^*=T$; 
 (f) the ITT augmented graph for $p(T^*,T, Y\mid F_T=t) = p(T^*,t, Y\mid F_T=t)$, 
  an intervention setting $T$ to $t$, so $F_T=t \neq \emptyset$;
(g) the latent projection of the graph in (d) after marginalizing $T$. 
Note that in (a), (b) we use $T^*$ (rather than $T$) for the natural value of treatment in order to highlight the correspondence to the ITT variables in Dawid's proposal. The graph in (g) corresponds to (a) and (b), under the correspondence $t \Leftrightarrow F_T\!=\!t$, $Y(t) \Leftrightarrow Y\,|\,F_T\!=\!t$.}
\label{fig:simple}}
\end{figure}
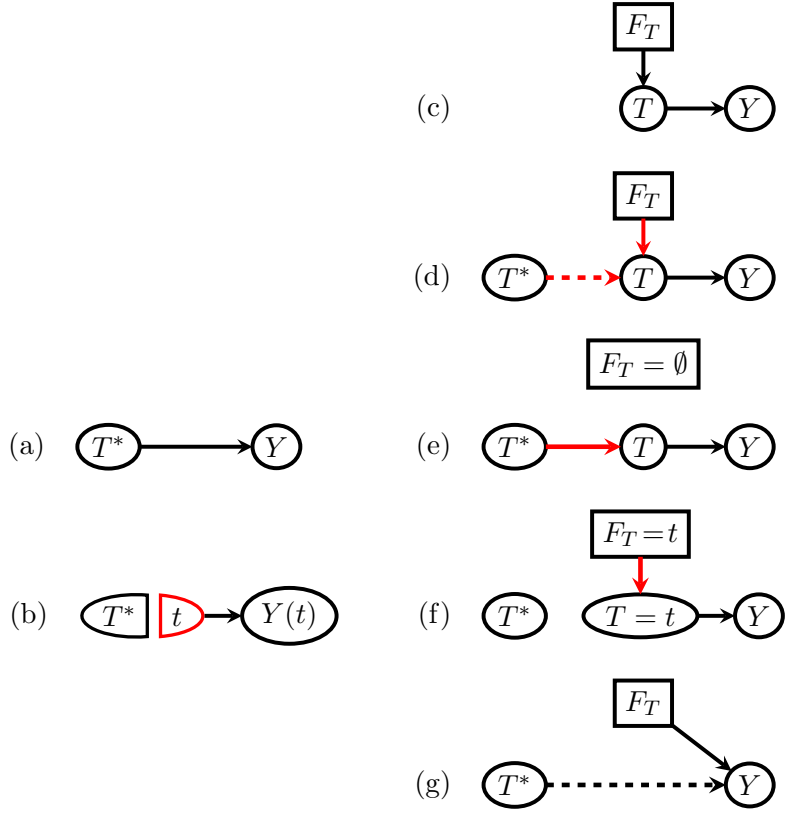

\eject

\section{Critique of Dawid's Proposal}

We have the following four main issues, which we describe in detail below:
\begin{itemize}
\item[(1)] The inclusion of ITT variables within Dawid's theory appears necessary in order to distinguish causal relationships from happenstance agreement between observational and (`fat hand') intervention distributions. However, including all three of  $T$ (the `actual' treatment), $T^*$ (the ITT variable),  and  $F_T$ (the regime indicator) introduces deterministically related variables and thereby obscures the content of Dawid's defining conditional independences A and B;
\item[(2)] Related to the previous point, d-separation is no longer a complete criterion for determining conditional independence on a graph in which there are definitional deterministic relationships between the variables;\footnote{This is also an issue for the Twin Network approach developed in \citet{pearl:2009}.}
\item[(3)] Dawid's ITT augmented diagrams incorporate context-specific independence (via dashed edges) but his results do not establish that the resulting distribution obeys all of the implied context-specific independences; these are not implied by his defining conditional independences $A+B$; 
these independences will not hold without additional information concerning the relation of  $T$ to $T^*$ and $F_T$ that is not captured in A+B;
\item[(4)] Dawid makes use of what he terms `fictitious' independence relations, but he argues that these are assumptions that can be made without loss of generality. This is not the case in general, though, as we show, in the context of his arguments the resulting logical `gap' can be filled.
\end{itemize}

We show that all of these issues may be avoided by re-formulating his theory in two simple ways:
\begin{itemize}
\item[(I)] Marginalizing out the post-intervention treatment variable $T$ while keeping the ITT variable $T^*$;\footnote{Dawid instead proposes
to marginalizes out the ITT variables.}
\item[(II)] Formulating the defining extended independence relations in terms of distributional consistency and the augmented ITT diagram (after marginalizing $T$) and intervening on all the variables in $A$; the local Markov property for the original variables is then implied.
\end{itemize}
The resulting theory is formally isomorphic to the SWIG theory described above; the augmented ITT graph can be viewed as containing the union
of the nodes and edges in the original DAG $\G$ and the SWIG $\G({a})$, with the fixed nodes in the SWIG corresponding to the (non-idle) regime indicators in the augmented DAG.

\subsection{The simplest setting}

Consider the setting in which there is a single exposure $T$ and an outcome $Y$; suppose that
$T$ takes a finite set of states $\frak{T}$.
Dawid's augmented causal graph with the Intention-To-Treat variable $T^*$ is shown in Figure \ref{fig:simple}(d). 
Here $T^*$ represents the natural value of treatment which an individual is ``selected to receive'' \citep[p.52]{dawid:2021} in the absence of an intervention that would override this. This is distinct from $T$ the ``treatment actually applied'' \citep[p.54, Def.~1]{dawid:2021}; 
$F_T$ is a regime indicator taking values in
$\frak{T} \cup \{\emptyset\}$. Under Dawid's proposal the graph in Figure \ref{fig:simple}(d) represents the kernel $p(T^*,T,Y \mid F_T)$;
$F_T=\emptyset$ indicates the observational regime in which case $T=T^*$, see Figure \ref{fig:simple}(e) where we have used a colored edge, $T^*\textcolor{red}{\rightarrow} T$, to indicate the deterministic relationship between $T$ and $T^*$. Similarly, 
$F_T=t\in \frak{T}$ indicates the interventional regime in which case $T=t$, see Figure \ref{fig:simple}(f). Note that $T\indep T^* \mid F_T \neq \emptyset$, which is represented by the dashed edge from $T^*$ to $T$ in Figure \ref{fig:simple} (d) and by the absence of the edge between $T^*$ and $T$ in Figure \ref{fig:simple}(e).

For comparison, Figure \ref{fig:simple}(a) and (b), respectively show representations of the observed distribution $p(T^*,Y)$ and the joint distribution $p(T^*,Y(t))$; as suggested by the graphical structures, there is a close correspondence between these approaches when ITT variables are included in the decision theory graph. In what follows we will show that in fact, the two theories can be shown to be isomorphic up to labeling of variables.

{
\begin{table}
\centering
\begin{tabular}{lcc}
& Pot.~Outcome & Decision Theoretic\\[5pt]
\hline\\[-5pt]
Graph for observed data &${\G}$ & ITT~DAG, $F_T=\emptyset$\\[4pt]
Graph representing intervention on $T$ &${\G}(t)$ & ITT DAG, $F_T=t$\\[10pt]
Observed distribution &$p(T^*, Y)$ & $p(T^*,Y \mid F_T = \emptyset)$\\[10pt]
 \begin{tabular}{@{}l@{}}Distribution resulting from setting\\
 $T=t$ directly after observing $T^*$\end{tabular} &$p(T^*,Y(t))$ & $p(T^*,Y \mid F_T = t)$\\

\end{tabular}
\caption{Correspondence between the potential outcome / SWIG approach and the decision theoretic approach.
Here, in the potential outcome approach we use $T^*$ (rather than $T$) to denote the natural value of treatment so as to make the correspondence more self-evident.
\label{tab:swig-decision}
}

\end{table}
}

Although Dawid includes ITT variables in the development here, they were absent in \citep{dawid:2000} and 
ultimately his goal is to remove the ITT variables, leaving the DAG shown in Figure \ref{fig:simple}(c) containing only the original variables and the 
treatment indicators; see bottom of \citep[p.65]{dawid:2021}. Dawid states that the augmented graphs without ITT variables are sufficient for reasoning about point interventions.

Given this, one may ask why it is necessary introduce the ITT variables into the theory in the first place.
One issue that arises is that without the ITT variables, the decision theoretic approach lacks the language to describe concepts such as the 
effect of treatment on the treated.  In addition, the approach lacks the concepts necessary to distinguish agreement between distributions in the observed and interventional world that reflect agreement between an observational study and a randomized experiment due to the absence of confounding, from an equality that is purely `contingent' or spurious. 

To illustrate this, consider the following story. Suppose that a manufacturer of dietary supplements carries out an observational study. They find that those who regularly consume the supplement ($T=1$) have lower levels of `bad' cholesterol ($Y$) than the people who do not ($T=0$). Buoyed by these results the manufacturer hires a company to perform a randomized trial. The results of the previous study are given to the company; it is made clear that the manufacturer would like these results confirmed and that repeat business depends on the firm achieving this.
In order to comply with this the testing company carries out a non-blinded study and also modify the software in the cholesterol measuring system to ensure that the results agree with those in the observational study; see Figure \ref{fig:coincidence}(b), here $H$ represents unobserved confounding and the edge $F_T \rightarrow Y$ indicates the compromised measurement process.\footnote{Within the potential outcome framework this would correspond to a failure of consistency, since among people with $T^*=t$, it need not hold that their observed outcome $Y$ is the same as the outcome they would have had, had they been in an experiment and assigned to $t$, namely $Y(t)$.} Since the experimental and observational distributions agree, it will hold that $Y \indep F_T \mid T$, as implied by the decision-theoretic graph in Figure  \ref{fig:coincidence}(a).

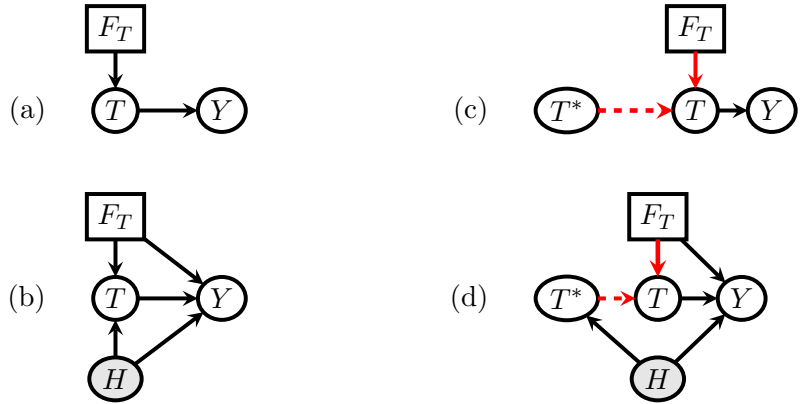
\begin{figure}
\centering
\begin{tikzpicture}
\tikzset{line width=1.25pt, outer sep=0pt,
         ell/.style={draw, inner sep=2pt,
          line width=1.5pt}, >=stealth,
         swig vsplit={gap=5pt,
         inner line width right=0pt, line color right=red}};
\begin{scope}
\node[name=tplus,ell,draw,shape=ellipse]{$T$};
\node[name=f,ell,draw,shape=rectangle, inner sep=4pt, above = 5mm of tplus]{$F_T$};
\node[name=y,ell,draw,shape=ellipse,right=8mm of tplus]{$Y$};
\draw[ell,->](tplus) to (y);
\draw[ell,->](f) to (tplus);
\node[left = 5mm of tplus]{(a)};
\end{scope}
\begin{scope}[yshift=-2.5cm]
\node[name=tplus,ell,draw,shape=ellipse]{$T$};
\node[name=f,ell,draw,shape=rectangle, inner sep=4pt, above = 5mm of tplus]{$F_T$};
\node[name=h,ell,draw,shape=ellipse, fill=black!10, inner sep=2pt, below = 5mm of tplus]{$H$};
\node[name=y,ell,draw,shape=ellipse,right=8mm of tplus]{$Y$};
\draw[ell,->](tplus) to (y);
\draw[ell,->](f) to (tplus);
\draw[ell,->](h) to (y);
\draw[ell,->](h) to (tplus);
\draw[ell,->](f) to (y);
\node[left = 5mm of tplus]{(b)};
\end{scope}
\begin{scope}[xshift=6cm]
\node[name=t,ell,draw,shape=ellipse]{$T^*$};
\node[name=tplus,ell,draw,shape=ellipse, right = 10mm of t]{$T$};
\node[name=f,ell,draw,shape=rectangle, inner sep=4pt, above = 5mm of tplus]{$F_T$};
\node[name=y,ell,draw,shape=ellipse,right=4mm of tplus]{$Y$};
\draw[ell,->,red, dashed, line width=1.75pt](t) to (tplus);
\draw[ell,->](tplus) to (y);
\draw[ell,red,->](f) to (tplus);
\node[left = 5mm of t]{(c)};
\end{scope}
\begin{scope}[xshift=6cm,yshift=-2.5cm]
\node[name=t,ell,draw,shape=ellipse]{$T^*$};
\node[name=tplus,ell,draw,shape=ellipse, right = 5mm of t]{$T$};
\node[name=f,ell,draw,shape=rectangle, inner sep=4pt, above = 5mm of tplus]{$F_T$};
\node[name=y,ell,draw,shape=ellipse,right=5mm of tplus]{$Y$};
\node[name=h,ell,draw,shape=ellipse, fill=black!10, inner sep=2pt, below = 5mm of tplus]{$H$};
\draw[ell,red, ->,line width=1.75pt](f) to (tplus);
\draw[ell,red, dashed,->](t) to (tplus);
\draw[ell,->](tplus) to (y);
\draw[ell,->](h) to (y);
\draw[ell,->](h) to (t);
\draw[ell,->](f) to (y);
\node[left = 5mm of t]{(d)};
\end{scope}
%
\end{tikzpicture}
\caption{Illustration of the necessity of the Intention-To-Treat (aka `natural value of treatment') variable $T^*$ in Dawid's proposal.
(a) An augmented DAG (without ITT nodes) corresponding to an observational study without confounding and a perfect intervention on $T$.
(b) An augmented DAG (without ITT nodes) representing an observational study with confounding ($H$) and a mis-targeted (`fat hand') intervention  affecting both $T$ and $Y$.
 If the mis-targeted intervention matches the effect of confounding then there will be equality of the observational and interventional distributions $p(Y|T=t, F_T=\emptyset) = p(Y|T=t, F_T=t)$ so that the extended independence $Y \indep F_T \mid T$ will hold and hence the causal diagram shown in (a) cannot be refuted. The inclusion of $T^*$ resolves this. (c) The DAG with ITT variables corresponding to the study without confounding, this implies $Y \indep F_T, T^* \mid T$, which is not implied by the ITT augmented DAG (d) when confounding is present.
 \label{fig:coincidence}
}
\end{figure}

To be clear, the critique here is {\it not} that someone who was unaware of the presence of confounding and the devious activities of the company running the trial would infer the wrong causal effect. Rather, it is that without the ITT variables, the decision-theoretic approach lacks the conceptual apparatus necessary to distinguish the situations in Figure \ref{fig:coincidence}(a) and (b).\footnote{Here we are assuming that there is no information available regarding the nature or identity of the possible confounding variables $H$.}
In contrast, if the ITT variables $T^*$ are included, then no such difficulty arises: the corresponding augmented DAG, shown in 
Figure  \ref{fig:coincidence}(c) now additionally requires that $Y \indep T^* \mid F_T=t$, which will fail to hold if there is unobserved confounding between $T^*$ and $Y$. Note that this latter condition is essentially equivalent to the ignorability condition $Y(t) \indep T^*$ in the potential outcome framework; we return to this point below.

\subsection{Dawid's defining extended conditional independence relations}

Under Dawid's formalism the augmented graph with ITT variables defines a causal model via the following extended conditional independence relations:
\begin{align}
A:&& T^*&\indep F_T \label{ind:a}\\
B:&& Y &\indep T^*, F_T \mid T \label{ind:b}
\end{align}
see \citet[Eq. (62), (63)]{dawid:2021}.

\subsubsection{Dawid's Independence A}\label{sec:dawid-a}
The first independence (\ref{ind:a}) states that whether or not there is an intervention on $T$ has no effect on the (distribution of the) ITT value $T^*$. 
Indeed, Dawid states:
\begin{quotation}
Now $T^*$ is determined prior to any (actual or hypothetical) treatment application, and behaves as a covariate [\ldots] this distribution is then the same in all regimes. \citep[\S8, p.54]{dawid:2021}.
\end{quotation}
Similarly, in the potential outcome framework, it is assumed that intervention on a treatment variable does not affect variables whose values are realized prior to that intervention, including the natural value of that treatment variable, $T^*$, so that $T^*(t) = T^*$.

However, Dawid's reference to $T^*$ being a covariate that is {\it determined prior} to an actual or {\it hypothetical} treatment application is perhaps surprising: If the value taken by $T^*$ is determined prior to the decision regarding the regime $F_T$ then this would appear to imply that in fact the random variables in the distributions $p(T^* \mid F_T=\emptyset)$ and $p(T^* \mid F_T=t)$ must live on a common probability space. But in this case it is hard to see why the random variables in the distributions 
$p(T^*,T,Y \mid F_T=\emptyset)$ and  $p(T^*,T,Y \mid F_T=t)$ should not also live on a common probability space! The primary obstacle to so doing appears to be the use of $Y$ and $T$ to indicate what are distinct random variables (corresponding to different regimes) that are defined on the same space. This problem can obviously be overcome by simply using $(T^*,T,Y)$ and $(T^*, T(t), Y(t))$ to refer to the random variables under the idle and intervention regimes respectively; following Definition 1 in \citep{dawid:2021}, this would imply that $T=T^*$ (under the idle regime) and $T(t) = t$ (under an intervention). 

An analyst who adopted this notation is not obligated to impose any additional equalities relating these random variables -- such as those implied by consistency -- should they not wish to do so. As we did above in \S\ref{sec:swigs}, one might choose instead to follow Dawid by merely imposing distributional consistency; see also further discussion below. 
However, from the perspective of the potential outcome framework this leads to an unnecessary multiplicity of random variables and more cumbersome notation.
For example, in the simple case of a binary treatment this approach requires three random variables $\{Y,Y(0),Y(1)\}$ corresponding to the response, rather than just two $\{Y(0),Y(1)\}$ with consistency at the level of random variables.\footnote{In other words, with $Y$ given by a deterministic function, so $Y=(1-T)Y(0)+TY(1)$.} It is unclear what is gained by assuming consistency at the level of distributions rather than individuals.

\subsubsection{Dawid's Independence B}
The fact that $T$ is a deterministic function of $T^*$ and $F_T$ means that the number of non-trivial conditional independence statements in (\ref{ind:b})  is not self-evident. A casual reader might imagine that in (\ref{ind:b}) the pair $(F_T, T^*)$ might take $(|\frak{T}|+1)|\frak{T}|$ different values for each value of the conditioning variable $T$. However, given $T=t$ there are only $|\frak{T}|+1$ possible values for $(F_T, T^*)$:
\[
T=t \quad \Rightarrow (F_T,T^*) \in \{(\emptyset, t)\} \cup \{ (t,s) : s \in \frak{T}\},
\]
since either we are in the idle regime, $F_T=\emptyset$ and $T=T^*$, or we are in the interventional regime, in which case $F_T=t$ and $T^*$ may take any value.
Thus given $T=t$, (\ref{ind:b}) corresponds to a set of $|\frak{T}|$ equalities:
\begin{align}
p(Y \mid T^*=t,F_T=t, T=t) &= p(Y \mid T^*=t,F_T=\emptyset, T=t),\label{eq:con1}\\
p(Y \mid T^*=t,F_T=t, T=t) &= p(Y \mid T^*=s,F_T=t, T=t),\quad \hbox{ for }s \neq t.\label{eq:ignor1}
\end{align}

\subsubsection*{Distributional Consistency in B} 
The equation (\ref{eq:con1}) corresponds to distributional consistency, which \citet[Eq.(14)]{dawid:2021} defines as:
\begin{align}
p(Y \mid F_T=\emptyset, T=t) &= p(Y \mid T^*=t,F_T=t).\label{eq:con4}
\end{align}
Dawid notes that this implies:
\begin{align}
Y&\indep F_T\mid T^*, T;\label{eq:dist-cons}
\end{align}
see \citet[Lemma 1]{dawid:2021}. However, this formulation also somewhat obscures the actual number of constraints: if $T^*\neq T=t$ then $F_T=t$ so that the statement becomes trivial, while if $T^*=T=t$ then $F_T$ only takes two possible values $\emptyset$ and $t$. Given this, it becomes clear that (\ref{eq:dist-cons}) may be reformulated by defining a dynamic regime $g^*$ that `intervenes' to set $T$ to be $T^*$. By defining a special regime indicator, denoted $F^*_T$, that takes only two values $\emptyset$ or $g^*$, we can re-express (\ref{eq:dist-cons}) as:
\begin{align}
Y, T^*&\indep F^*_T.\label{eq:dist-cons2}
\end{align}
Note that in so doing we do not need to refer to $T$;\footnote{Along similar lines, in his equation (14) \citet{dawid:2021} notes that the LHS of
(\ref{eq:con4}) is equivalent to $p(Y \mid F_T=\emptyset, T^*=t)$.}  see Figure \ref{fig:consistency-phil}(d) for a graphical depiction.

In terms of potential outcomes the independence (\ref{eq:dist-cons2}) may be expressed as:
 \begin{align}
p(Y, T^*=t) &= p(Y, T^*=t \mid F^*_T=\emptyset)\nonumber\\
&= p(Y, T^*=t \mid F^*_T=g^*)\nonumber\\
&= p(Y(t), T^*(t)=t),\label{eq:dist-cons3}
\end{align}
which corresponds to distributional consistency; see (\ref{eq:consist-swig-simple}).

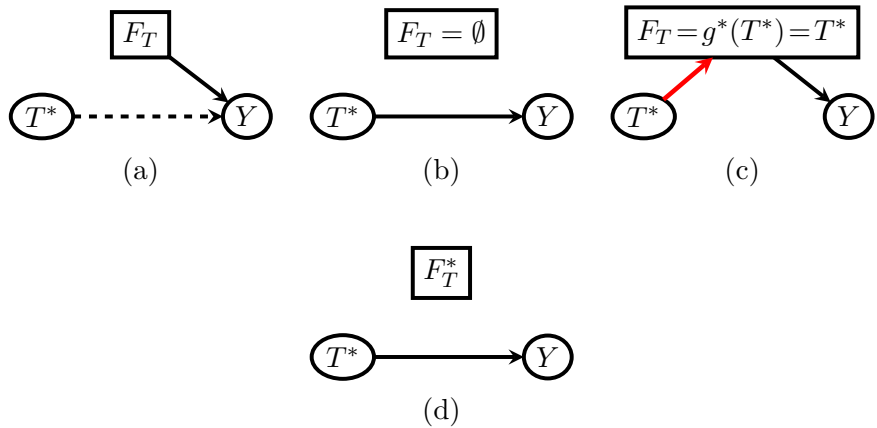
\begin{figure}
\centering
\begin{tikzpicture}
\tikzset{line width=1.25pt, outer sep=0pt,
         ell/.style={draw, inner sep=2pt,
          line width=1.5pt}, >=stealth,
         swig vsplit={gap=5pt,
         inner line width right=0pt, line color right=red}};
\begin{scope}
\node[name=t,ell,draw,shape=ellipse]{$T^*$};
\node[name=tplus,ell,draw=none,opacity=0,shape=ellipse, right = 6mm of t]{\phantom{$T$}};
\node[name=f,ell,draw,shape=rectangle, inner sep=4pt, above = 5mm of tplus]{$F_T$};
\node[name=y,ell,draw,shape=ellipse,right=8mm of tplus]{$Y$};
\draw[ell,->,dashed](t) to (y);
\draw[ell,->](f) to (y);
\node[below=12mm of f]{(a)};
\end{scope}
\begin{scope}[xshift=4cm,yshift=0cm]
\node[name=t,ell,draw,shape=ellipse]{$T^*$};
\node[name=tplus,ell,draw=none,opacity=0,shape=ellipse, right = 6mm of t]{\phantom{$T$}};
\node[name=f,ell,draw,shape=rectangle, inner sep=4pt, above = 5mm of tplus]{$F_T=\emptyset$};
\node[name=y,ell,draw,shape=ellipse,right=8mm of tplus]{$Y$};
\draw[ell,->](t) to (y);
\node[below=12mm of f]{(b)};
\end{scope}
\begin{scope}[xshift=8cm,yshift=0cm]
\node[name=t,ell,draw,shape=ellipse]{$T^*$};
\node[name=tplus,ell,draw=none,opacity=0,shape=ellipse, right = 6mm of t]{\phantom{$T$}};
\node[name=f,ell,draw,shape=rectangle, inner sep=4pt, above = 5mm of tplus]{$F_T\!=\!g^*(T^*)\!=\!T^*$};
\node[name=y,ell,draw,shape=ellipse,right=8mm of tplus]{$Y$};
\draw[ell,->,red, line width=1.75pt](t) to (f);
\draw[ell,->](f) to (y);
\node[below=12mm of f]{(c)};
\end{scope}
\begin{scope}[xshift=4cm,yshift=-3.2cm]
\node[name=t,ell,draw,shape=ellipse]{$T^*$};
\node[name=tplus,ell,draw=none,opacity=0,shape=ellipse, right = 6mm of t]{\phantom{$T$}};
\node[name=f,ell,draw,shape=rectangle, inner sep=4pt, above = 5mm of tplus]{$F^*_T$};
\node[name=y,ell,draw,shape=ellipse,right=8mm of tplus]{$Y$};
\draw[ell,->](t) to (y);
\node[below=12mm of f]{(d)};
\end{scope}
\end{tikzpicture}
\caption{Encoding Distributional Consistency via a special dynamic regime in the setting of a reformulated decision diagram (having marginalized the intervention target `$T$').
(a) Reformulated augmented graph ${\G}^*$ representing the observed joint distribution $p(T^*,Y \mid F_T)$; 
(b) Augmented graph ${\G}^*$ corresponding to $p(T^*,Y \mid F_T=\emptyset)$;
(c) Augmented graph ${\G}^*$ corresponding to $p(T^*,Y \mid F_T=g^*)$; corresponding to the dynamic `regime' that `intervenes' on the target  setting it to 
the natural value $T^*$; (d) A graph illustrating distributional consistency (\ref{eq:dist-cons2}); here $F^*_T$ is a special regime indicator taking only the values $\emptyset$ and $g^*$. The graph (d) encodes the distributional consistency assumption: the distribution over $Y$ and $T^*$ resulting from the `intervention' $g^*$ is identical to having no intervention.
\label{fig:consistency-phil}}
\end{figure}

\begin{figure}
\centering
\begin{tikzpicture}
\tikzset{line width=1.25pt, outer sep=0pt,
         ell/.style={draw, 
         inner sep=2pt,
          line width=1.5pt}, >=stealth,
         swig vsplit={gap=5pt,
         inner line width right=0pt, line color right=red}};
\begin{scope}
\node[name=tstar,ell,draw,shape=ellipse]{$T^*$};
\node[name=tplus,ell,draw=none,opacity=0,shape=ellipse, right = 4mm of tstar]{\phantom{$T$}};
\node[name=f,ell,draw,shape=rectangle, inner sep=4pt, above = 8mm of tstar]{$F_T$};
\node[name=y,ell,draw,shape=ellipse,right=14mm of tstar]{$Y$};
\draw[ell,->,dashed, line width=1.75pt](tstar) to (y);
\draw[ell,->](f) to[out=0,in=90] (y);
\node[below=5mm of tplus]{(a)};
\end{scope}
\begin{scope}[xshift=6cm]
\node[name=tstar,ell,draw,shape=ellipse]{$T^*$};
\node[name=tplus,ell,draw=none,opacity=0,shape=ellipse, right = 2mm of tstar]{\phantom{$T$}};
\node[name=f,ell,draw,shape=rectangle, inner sep=4pt, above = 8mm of tstar]{$F_T$};
\node[name=y,ell,draw,shape=ellipse,right=14mm of tstar]{$Y$};
\draw[ell,->](f) to[out=0,in=90] (y);
\node[name=t,ell,draw,shape=ellipse,above=1mm of tplus, xshift=2mm]{$T$};
\draw[ell,->, dashed, line width=1.75pt](tstar) to (y);
\draw[ell,red,->](f) to (t);
\draw[ell, dashed, red,->](tstar) to (t);
\node[below=5mm of tplus]{(b)};
\end{scope}
\end{tikzpicture}
\caption{ (a) Reformulated augmented graph ${\G}^*$ representing the observed joint distribution $p(T^*,Y \mid F_T)$; 
(b) graph illustrating that, if desired, the `applied treatment' variable $T$ may be added to $\G^*$ since it is a deterministic function of $T^*$ and $F_T$.
Note that although it may seem counterintuitive that $T$ is not a parent of $Y$ in this graph, this is formally correct. \label{fig:phil-with-t}}
\end{figure}
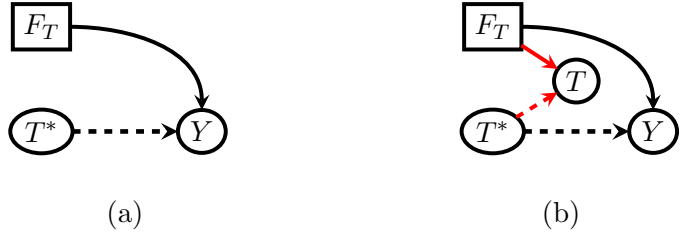

\subsubsection*{Ignorability in B} 
The equations (\ref{eq:ignor1}) express the property of ignorability, which \citet[eq.(20)]{dawid:2021} expresses as:
\begin{align}
Y&\indep T^*\mid F_T, T.\label{eq:dist-ignorability}
\end{align}
However, as Dawid himself notes, given $T=t$, then either $F_T=\emptyset$ in which case $T^*=t$ (and independence holds trivially),
or $F_T=t$, so that this constraint is identical to:
\begin{align}
Y&\indep T^*\mid F_T=t\quad \hbox{ for } t\in \frak{T};\label{eq:dist-ignorability2}
\end{align}
see Figure \ref{fig:simple}(f),(g).
Equivalently in terms of potential outcomes,
\begin{align}
Y(t) &\indep T^* \hbox{ for } t\in \frak{T};\label{eq:dist-ignorability3}
\end{align}
see Figure \ref{fig:simple}(b).
Again, we note that  $T$ is not required for the purpose of expressing this condition.

\subsection{Simplification}
From a graphical perspective it is perhaps natural to wish to express the invariance of the distribution of $Y$ given $T$ across observational and interventional distributions by examining whether a regime indicator $F_T$ is d-separated from $Y$ given $T$. However, as we have seen, it is necessary to include what Dawid calls the ITT variable (aka the natural value of treatment) $T^*$ in order to rule out cases of spurious invariance.
Furthermore, $T^*$ plays a central role in certain notions, such as the effect of treatment on the treated, that are widely used in many studies that apply the potential outcome framework.

As shown above, there is no need to condition on $T$ when describing the defining independences, and in fact doing so arguably obscures the  nature of the specific assumption being made. This suggests that $T$ should be marginalized from the ITT augmented graph, rather than  $T^*$ as Dawid proposes. Note that, if we distinguish the cases $F_T=\emptyset$ and $F_T=t$, the resulting graphs (modulo labeling) are isomorphic to those used in the SWIG framework; compare Figure \ref{fig:simple}(a) to (e), and (b) to (f). 

We carry out this reformulation in full generality in the next section.

 \section{Reformulation of Decision Graphs}
 
Our proposed reformulation of Decision Graphs follows a strategy similar to that used for SWIGs.
In contrast, Dawid aims to give extended conditional independence relations that, together with the usual independence relations over the observed variables, will yield the Markov property for the augmented Decision Graph with ITT variables. As an alternative, we begin by defining a Markov property associated with the augmented decision graph, and then, using distributional consistency derive the usual observed conditional independences.\footnote{However, see Appendix \ref{app:dawid-swig} for an alternative reformulation that, similar to Dawid's approach, starts from the assumption that the observed distribution $p(V)$ obeys the usual Markov property for $\G$.}
 
It should be noted that Dawid's independences do not actually imply the full Markov property for the ITT graph because, as noted by the presence of a dashed edge, there are context specific independences implied by the graph.\footnote{Note that in Figure 12 in \citep{dawid:2021},
there is an edge $Z\rightarrow X_1$ that is not dashed, but which it appears should be dashed: note that the relationship of $Z$ and $H$ to $X_1$ are symmetric and the edge $H \rightarrow X_1$ is dashed.}
 However, these are not implied by the independence relations $A$ and $B$. (To see this, note that the conditions $A$ and $B$ would also hold for a decision DAG with the same structure, but in which $T$ was not a deterministic function of $T^*$ and $F_T$, in which case the context specific independence relations would not hold.)
 
Note also that these extra extended conditional independence relations are not restricted solely to those involving $T$. Consider, for example, the 
front door graph shown in Figure  \ref{fig:front-door}(a). 
Since Dawid's augmented decision diagram in Figure  \ref{fig:front-door}(c) includes a dashed edge from $T^*$ to $T$, indicating that this edge should be removed conditional on $F_T=t$, the diagram implies that $Y$ will be d-separated from $F_T$ given $M$ and $F_T\neq \emptyset$. However, even though it is encoded in the augmented graph, the corresponding extended conditional independence:
\begin{align}\label{eq:front-door-specific}
Y \indep F_T \mid M, F_T \neq \emptyset,
\end{align}
does not follow from the independences A+B.

In the potential outcome framework, the constraint (\ref{eq:front-door-specific}) corresponds to:
\begin{align}\label{eq:swig-indep-front-door}
p(Y(t) \mid M(t)) = p(Y(t^*) \mid M(t^*)).
\end{align}
This constraint is naturally encoded by the d-separation of $Y(t)$ from the  fixed variable $t$ given $M(t)$ on the SWIG $\G(t)$ shown in Figure \ref{fig:front-door}(b); see \cite{shpitser2021multivariate, malinsky19po,robins2018dseparation}.

As these examples suggest, in order to capture the full Markov structure of the augmented decision diagram, including those constraints corresponding to dashed edges, 
it is natural to use the constraints implied by the decision diagram when no regime indicators are idle, which we express in shorthand as $F_A\neq \emptyset$; graphically this corresponds to removing (temporarily) all of the dashed edges. We show below that the independences encoded then imply,  via distributional consistency, the Markov property for the observed data that is encoded in the original graph.

Another advantage of this approach is that we will only require the ITT variables $T^*$; the `applied treatment' which \citet{dawid:2021} denotes `$T$' will not be required.%
\footnote{Since $T$ is a deterministic function of $T^*$ and $F_T$, it is possible to add back in these variables if we wish; see Figure \ref{fig:phil-with-t}.}

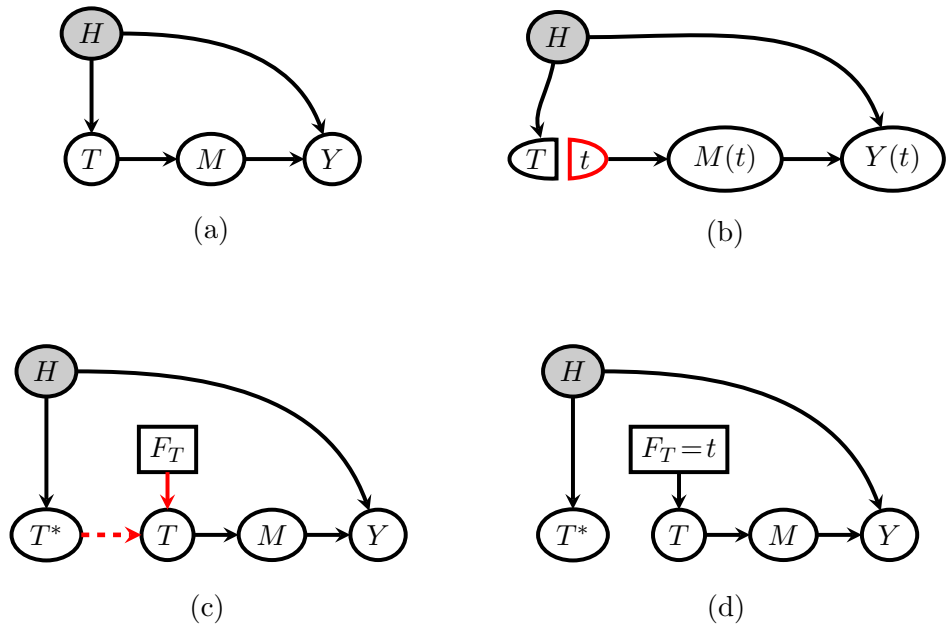
\begin{figure}
\begin{center}
\begin{tikzpicture}[
pre/.style={->,>=stealth,semithick,blue,line width = 1pt}]
\tikzset{line width=1.25pt, outer sep=0pt,
         ell/.style={draw, inner sep=3pt,
          line width=1.5pt}, >=stealth,
         swig vsplit={gap=5pt,
         inner line width right=0pt, line color right=red},
         ell2/.style={draw,fill opacity=0, text opacity=1}};
\begin{scope}[xshift=0.6cm, yshift=5cm]
\node[name=tplus,ell,draw,shape=ellipse]{$T$};
\node[name=h,ell,draw,fill=black!20, shape=ellipse, above = 10mm of tplus]{$H$};
\node[name=m,ell,draw,shape=ellipse,right=8mm of tplus]{$M$};
\node[name=y,ell,draw,shape=ellipse,right=8mm of m]{$Y$};
\draw[ell,->](tplus) to (m);
\draw[ell,->](h) to (tplus);
\draw[ell,->](h) to[out=0,in=110](y);
\draw[ell,->](m) to (y);
\node[below=3mm of m]{(a)};
\end{scope}
\begin{scope}[xshift=6.8cm, yshift=5cm]
\node[name=tplus, ell,shape=swig vsplit, ell2]{
        \nodepart{left}{$T$}
        \nodepart{right}{$t$}};
\node[name=h,ell,draw,fill=black!20, shape=ellipse, above = 10mm of tplus]{$H$};
\node[name=m,ell,draw,shape=ellipse,right=8mm of tplus, inner sep=3pt]{$M(t)$};
\node[name=y,ell,draw,shape=ellipse,right=8mm of m, inner sep=3pt]{$Y(t)$};
\draw[ell,->](tplus) to (m);
\draw[ell,->](h) to[out=260,in=110] (tplus);
\draw[ell,->](h) to[out=0,in=110](y);
\draw[ell,->](m) to (y);
\node[below=2.5mm of m]{(b)};
\end{scope}
\begin{scope}
\node[name=t,ell,draw,shape=ellipse]{$T^*$};
\node[name=tplus,ell,draw,shape=ellipse, right = 8mm of t]{$T$};
\node[name=f,ell,draw,shape=rectangle, inner sep=4pt, above = 5mm of tplus]{$F_T$};
\node[name=h,ell,draw,fill=black!20, shape=ellipse, above = 15mm of t]{$H$};
\node[name=m,ell,draw,shape=ellipse,right=6mm of tplus]{$M$};
\node[name=mplus,shape=ellipse,right=0mm of tplus]{\vphantom{$M$}};
\node[name=y,ell,draw,shape=ellipse,right=6mm of m]{$Y$};
\draw[ell,red,->,dashed, line width=1.75pt](t) to (tplus);
\draw[ell,->](tplus) to (m);
\draw[ell,red,->](f) to (tplus);
\draw[ell,->](h) to (t);
\draw[ell,->](h) to[out=0,in=110](y);
\draw[ell,->](m) to (y);
\node[below=3mm of mplus]{(c)};
\end{scope}
\begin{scope}[xshift=7cm, yshift=0cm]
\node[name=t,ell,draw,shape=ellipse]{$T^*$};
\node[name=tplus,ell,draw,shape=ellipse, right = 6mm of t]{$T$};
\node[name=f,ell,draw,shape=rectangle, inner sep=4pt, above = 5mm of tplus]{$F_T\!=\!t$};
\node[name=h,ell,draw,fill=black!20, shape=ellipse, above = 15mm of t]{$H$};
\node[name=m,ell,draw,shape=ellipse,right=6mm of tplus]{$M$};
\node[name=mplus,shape=ellipse,right=1mm of tplus]{\vphantom{$M$}};
\node[name=y,ell,draw,shape=ellipse,right=6mm of m]{$Y$};
\draw[ell,->](tplus) to (m);
\draw[ell,->](f) to (tplus);
\draw[ell,->](h) to (t);
\draw[ell,->](h) to[out=0,in=110](y);
\draw[ell,->](m) to (y);
\node[below=3mm of mplus]{(d)};
\end{scope}
\end{tikzpicture}
\end{center}
\caption{(a) Front-door graph $\G$; (b) the SWIG $\G(t)$ (with ancestral labeling); (c) the augmented decision diagram $\G^*$; (d) the augmented decision diagram given $F_T=t$ in which
the dashed edge from $T^*$ to $T$ is removed.
Notice that in (d) $F_T$ is d-separated from $Y$ given $M$. However, the corresponding extended independence, $Y \indep F_T \mid M, F_T\neq \emptyset$  is not implied by Dawid's conditions A+B.
\label{fig:front-door}}
\end{figure}


{\begin{table}
\centering
  \begin{tabular}{r@{\extracolsep{2pt}}c@{\extracolsep{2pt}}l}
 \multicolumn{3}{c}{Reformulated Decision Diagram Local Property}\\[5pt]
 \hline\\[-10pt]
 $H \indep$ & ${F_0},{F_1}$& $ \mid F_{01}\!\neq\!\emptyset$\\[2pt]
 $X_0 \indep$ & $ {H}, {F_0},{F_1}$& $ \mid F_{01}\!\neq\!\emptyset$\\[2pt]
 $Z \indep$ & $ {X_0}, {F_1}$& $ \mid H, F_0, F_{01}\!\neq\!\emptyset$\\[2pt]
  $X_1 \indep$ & $ {X_0}, {F_0}, {F_1} $& $ \mid H, Z, F_{01}\!\neq\!\emptyset$\\[2pt]
  $Y \indep$ & $ {H}, {X_0},  {X_1}, {F_0}$& $ \mid  Z, F_1, F_{01}\!\neq\!\emptyset$
 \end{tabular}
\caption{Defining properties for reformulated Decision Diagram corresponding to Figure \ref{fig:pearl-wrong-example2}, under the ordering $(H,X_0,Z,X_1,Y)$.
Here, as elsewhere, $F_{01}\neq \emptyset$ is a shorthand for $(F_0\neq \emptyset\,\, \&\,\, F_1\neq \emptyset)$.
\label{tab:pearl-wrong-example2}}
\end{table}
}

Specifically, consider a set of variables $V_1,\ldots ,V_p$.
 Let $A\subset \{1,\ldots, p\}$ be the (index) set of the targets of intervention.
If $i \in A$ then let $V_i$ be the corresponding ITT variable, (which Dawid denotes by $X^*_i$).
Thus the set $V_1,\ldots ,V_p$ consists of ITT variables as well as variables that are not in $A$ and hence not targets of intervention.\footnote{In this formulation we will not use the post-intervention target variables, which Dawid denotes, $X_i$.} 
Thus under the regime where every intervention target has been intervened upon, so that  $F_{i}\neq\emptyset$ for all $i \in A$,
the variables in $V_1,\ldots , V_p$  correspond to the random variables in the SWIG $\G(a)$.

For every intervention target $i\in A$, let $g_i^*$ denote the dynamic regime which `intervenes' to set the intervention target to its natural value $V_i$. Let $F^*_i$ be a regime indicator taking the states $\emptyset$ or $g_i^*$.

\begin{definition}[Distributional Consistency for Decision Diagrams]\label{def:dist-cons-dd}
The kernel $p(V\,|\, F_A)$ is said to obey distributional consistency if, 
given $B_i\in A$ and $C\subseteq A\setminus \{B_i\}$, where $C$ may be empty, 
\begin{align}\label{eq:consist-ind}
V \indep F^*_i \mid F_C, F_{C}\neq \emptyset,
\end{align}
where we use the shorthand
$F_C\neq \emptyset$ to indicate that for all $j \in C$, $F_j \neq \emptyset$.
\end{definition}

\noindent Note that, taking $Y=V\setminus \{B_i\}$,  (\ref{eq:consist-ind}) is equivalent to the following equality, which corresponds exactly to (\ref{eq:consist-swig}):
\begin{align}\label{eq:consist-ind2}
\MoveEqLeft{p(Y=y, B_i = b \mid F_{i} = b, F_{C}=c)}\\
&= p(Y=y, B_i = b \mid F^*_{i} = g^*_i, F_{C}=c, F_{C}\neq \emptyset)\nonumber \\
&=  p(Y=y, B_i = b \mid  F^*_{i} = \emptyset, F_{C}=c, F_{C}\neq \emptyset)\nonumber\\
&=  p(Y=y, B_i = b \mid  F_{i} = \emptyset, F_{C}=c)\nonumber\\
&=  p(Y=y, B_i=b \mid F_{C}=c).\label{eq:consist-ind2b}
\end{align}
Here the second equality follows from (\ref{eq:consist-ind}) while the first and third equalities are via the definition of $F^*$ and $\emptyset$.

As observed by Dawid, in place of Definition \ref{def:dist-cons-dd} we could instead have defined Distributional Consistency, without reference to the dynamic regime $g_i^*$, by simply equating (\ref{eq:consist-ind2}) and (\ref{eq:consist-ind2b}). We have chosen to make use of $g_i^*$ in order to emphasize what we see as the tautological nature of Distributional Consistency, while also formulating the condition (\ref{eq:consist-ind}) as a conditional independence.

The following three Lemmas are re-formulations of Lemmas \ref{lem:dist-consistency-vector}--\ref{lem:remove-b-from-conditional-not-a-fn} in the decision diagram framework. Though the proofs are largely translations of those lemmas, we include them here for completeness.

\begin{lemma}\label{lem:dist-consistency-vector-phil}
If $p(V\,|\, F_A)$ obeys distributional consistency,   $B$, $C$ are disjoint subsets of $A$, where $C$ may be empty,
then: 
\begin{align}\label{eq:consist-vector-phil}
V \indep F^*_B \mid F_C, F_{C}\neq \emptyset,
\end{align}
where $F^*_B$ is the set $\{F^*_i, i \in B\}$.
\end{lemma}

\noindent{\it Proof:} This follows by induction on the size of $B$.\hfill$\Box$

\begin{lemma}\label{lem:cond-consistency-phil} Let $B$, $C$ be disjoint subsets of $A$, where $C$ may be empty, and let $Y$, $W$ be disjoint subsets of $V\setminus B$, then distributional consistency implies:
\begin{align}\label{eq:consist-cond-swig-phil}
Y \indep F^*_B \mid B, W, F_C, F_C\neq \emptyset.
\end{align}
\end{lemma}

\noindent{\it Proof:} This follows by applying (extended) graphoid axioms to (\ref{eq:consist-vector-phil}).\hfill$\Box$

\begin{lemma}\label{lem:remove-b-from-joint-not-a-fn-phil}
Let $B$, $C$ be disjoint subsets of $A$, where $C$ may be empty. If $B\subseteq W$, then under distributional consistency: 
\begin{equation}\label{eq:remove-b-from-joint-not-a-fn-phil}
W \indep F_B \mid F_C, F_{B\cup C}\neq \emptyset \quad \Rightarrow\quad W \indep F_B \mid F_C, F_C\neq\emptyset.
\end{equation}
\end{lemma}

\noindent{\it Proof:} Let $b\in \mathfrak{X}_B$, $c\in \mathfrak{X}_C$ and $w\in \mathfrak{X}_W$. 
 Given the LHS of (\ref{eq:remove-b-from-joint-not-a-fn-phil})
it is sufficient to prove  $p(X_{W}\mid F_B=b,F_C=c) = p(X_{W}\mid F_B=\emptyset,F_C=c)$. Now:
\begin{align*}
\MoveEqLeft{p(X_{W}=w\mid F_B=b,F_C=c)}\\
 &= p(X_{W\setminus B}=w_{W\setminus B}, X_{B}=w_{B} \mid F_B=b,F_C=c)\\ 
&= p(X_{W\setminus B}=w_{W\setminus B}, X_{B}=w_{B}\mid F_B=w_B, F_C=c)\\
&= p(X_{W\setminus B}=w_{W\setminus B}, X_{B}=w_{B}\mid F^*_B=g^*_B,F_C=c)\\
&= p(X_{W\setminus B}=w_{W\setminus B}, X_{B}=w_{B}\mid F^*_B=\emptyset, F_C=c)\\
&= p(X_{W\setminus B}=w_{W\setminus B}, X_{B}=w_{B}\mid F_B=\emptyset, F_C=c).
\end{align*}
The second equality uses the premise of (\ref{eq:remove-b-from-joint-not-a-fn-phil}), the third is by definition of $g^*_B$, the fourth
is distributional consistency via Lemma \ref{lem:dist-consistency-vector-phil}. \hfill$\Box$

\begin{lemma}\label{lem:remove-b-from-conditional-not-a-fn-phil}
Let $B$, $C$ be disjoint subsets of $A$, where $C$ may be empty.  Let $Y$, $W$ be disjoint sets with $B\subseteq W$,
then under distributional consistency:
\begin{equation}\label{eq:remove-b-from-conditional-not-a-fn-phil}
Y \indep F_B \mid W, F_C, F_{B\cup C}\neq \emptyset \quad \Rightarrow\quad Y \indep F_B \mid W, F_C, F_C\neq\emptyset.
\end{equation}
\end{lemma}

\noindent{\it Proof:} Similar to the proof of Lemma \ref{lem:remove-b-from-joint-not-a-fn-phil}, 
given the premise in (\ref{eq:remove-b-from-conditional-not-a-fn-phil}), it suffices to show that
$p(X_{Y}\mid X_W, F_B=b, F_C=c) = p(X_{Y}\mid X_W, F_B=\emptyset,F_C=c)$.
\begin{align*}
\MoveEqLeft{p(X_{Y} \mid X_W=w, F_B=b, F_C=c)}\\
 &= p(X_{Y} \mid X_{W\setminus B}=w_{W\setminus B}, X_{B}=w_{B}, F_B=b,F_C=c)\\ 
&= p(X_{Y} \mid X_{W\setminus B}=w_{W\setminus B}, X_{B}=w_{B}, F_B=w_{B}, F_C=c)\\ 
&= p(X_{Y} \mid X_{W\setminus B}=w_{W\setminus B}, X_{B}=w_{B}, F^*_B=g^*_B, F_C=c)\\ 
&= p(X_{Y} \mid X_{W\setminus B}=w_{W\setminus B}, X_{B}=w_{B}, F_B=\emptyset, F_C=c).
\end{align*}
As in the previous proof, the second equality uses the premise of (\ref{eq:remove-b-from-conditional-not-a-fn-phil}), the third is by definition of $g^*_B$, the fourth is distributional consistency via Lemma \ref{lem:dist-consistency-vector-phil}. \hfill$\Box$
\\

\subsection{Reformulated Augmented Decision Diagrams}

Let $\G$ be a DAG with a topologically ordered vertex set $V=\{1,\ldots ,p\}$  representing an observed distribution $p(W_V)$.\footnote{In this section of the paper, when we wish to distinguish random variables from index sets, we use $W_i$, rather than $X_i$. This is because
in Dawid's development $X$ is reserved to denote intervention targets. However, we depart from Dawid's notation in that we will not use an asterisk to indicate ITT variables, this is because we will only include the ITT variables associated with intervention targets on the reformulated diagram.}
Let $A \subseteq V$ be the subset of vertices for which interventions are well-defined, let $F = \{F_i, i \in A\}$ be the corresponding set of regime indicators.  Let $\G^*$ be the extended DAG with vertex set $V\cup F$, representing the kernels $p(W_V \mid F_A)$.
As before we use $\pa(i)$ to indicate the (index) set of variables that are parents of $W_i$ in the original DAG $\G$, and let
$\pre(i)$  denote  $\{1,\ldots ,i-1\}$, the predecessors of $i$ under a total ordering $\prec$ consistent with $\G$.\footnote{Note that the vertex set for $\G^*$ corresponds to the sets of variables that in Dawid's notation would be written $(V_i,X_1^*,\ldots V_k, X_k^*,V_{k+1})$; in other words, it consists of domain variables and ITT variables associated with intervention targets. We will have no need to include what Dawid calls `the intervention targets' which he denotes $X_i$, in our reformulated decision diagram, though they may be added; see Figure \protect\ref{fig:phil-with-t}.}

\begin{definition}\label{def:local-swig-mp-phil}
The kernel $p(W_V \mid F_A)$ will be said to obey the augmented DAG local Markov property for the DAG $\G^*$ if for all $i\in V$:
\begin{align}\label{eq:local-ind}
W_i \indep F_{A\setminus \pa(i)}, W_{\pre(i) \setminus (\pa(i)\setminus A)}
\mid W_{\pa(i) \setminus A}, F_{A \cap \pa(i)},  F_{A} \neq \emptyset,
\end{align}
where $F_{A\cap \pa(i)} \neq \emptyset$ is a shorthand for $F_j \neq \emptyset$ for all $j \in A\cap \pa(i)$.
\end{definition}

This formulation captures the Markov property necessary for the augmented diagram including the context specific independences that arise from interventions (that are not captured directly in Dawid's A+B formulation).

Note that this property follows from d-separation applied to the graph in which we intervene on every vertex in $A$.
 We will show that under distributional consistency this property 
implies factorization of the observed distribution with respect to the original graph.

However, it is useful first to further decompose the sets on the RHS of the independence. Specifically, we divide the regime indicators that are not parents of $i$ into those that occur after $i$ and those that are prior to $i$:
\[
F_{A\setminus \pa(i)} = \left( F_{A\setminus \pre(i)}, F_{(A\cap \pre(i))\setminus \pa(i)}\right).
\]
Similarly, we divide the set of random variables that are prior to $i$ and either in $A$ or not parents of $i$ into 
those that are not parents and those that are parents that are in $A$: 
\[
W_{\pre(i)\setminus (\pa (i)\setminus A)} =  \left( W_{\pre(i)\setminus \pa(i)}, W_{\pa (i)\cap A}\right).
\]
Thus independence (\ref{eq:local-ind}) becomes:
\begin{align}
W_i &\indep \overbrace{F_{A\setminus \pre(i)}}^{\hbox{\tiny time order}},\;\;   
\overbrace{F_{(A\cap \pre(i))\setminus \pa(i)}}^{\hbox{\tiny causal Markov prop.}},\;\;\overbrace{W_{\pre(i)\setminus \pa(i)},}^{\hbox{\tiny assoc.~Markov prop.}} \quad \overbrace{W_{\pa(i) \cap A}}^{\hbox{\tiny ignorability}}\label{eq:local-ind2} \\[4pt]
&\kern120pt  
\left| \vphantom{\bigcup}\right.
\underbrace{F_{A \cap \pa(i)},}_{\hbox{\tiny~~ fixed parents~~}}
\underbrace{W_{\pa(i) \setminus A},}_{\hbox{\tiny ~~random parents~~}}
\underbrace{F_{A\vphantom{\mid}} \neq \emptyset.}_{\hbox{\tiny ~~ intervene on all of $A$~~}}\nonumber
\end{align}

Consequently, independence (\ref{eq:local-ind2}) captures the following:
\begin{itemize}
\item later interventions have no effect on earlier distributions (time order);
\item given intervention on all earlier targets, the specific value of an intervention does not affect the distribution of a variable given its non-intervened parents unless the intervened on variable is itself a parent (causal Markov property);
\item independence from earlier random variables given non-intervened parents (associational Markov property);
\item an intervention on a parent of a variable renders that variable independent of the natural value of the intervention target conditional on its other non-intervened parents (ignorability).
\end{itemize}

\subsection{Example}

In Table \ref{tab:pearl-wrong-example2} we show the reformulated Decision Diagram Markov property corresponding to the augmented DAG $\G^*$ shown in Figure \ref{fig:pearl-wrong-example2}(c). Notice that the local property here corresponds naturally to the graph $\G^*$ under the regime $F_0=x_0$, $F_1=x_1$ displayed in Figure \ref{fig:pearl-wrong-example2}(d).
In particular, note that for each random vertex, the size of the conditioning set in the defining independence (ignoring the term $F_{01}\neq \emptyset$) is equal to the number  of parents that the vertex has in Figure \ref{fig:pearl-wrong-example2}(d).

\subsection{Consequences of the Local Markov property}

\begin{lemma}\label{lem:kernel-of-swig-lmp-phil}
If the kernel $p(W_V \mid F_A)$ obeys distribution consistency and the augmented DAG local Markov property w.r.t. $\G^*$ then:
\begin{align}
\MoveEqLeft{p(W_i \mid W_{\pre(i)}, F_A\!=\!a )}\label{eq:swig-lmp-lhs-phil}\\
&= p(W_i \mid W_{\pre(i)}, F_{\pre(i)\cap A}\!=\!a_{\pre(i)\cap A} ) \label{eq:drop-successors-from-a-phil}  \\
&= p(W_i \mid W_{\pre(i)}, F_{\pa(i)\cap A}\!=\!a_{\pa(i)\cap A} )\label{eq:drop-predecessors-from-a-phil} \\
&= p(W_i \mid W_{\pa(i)}, F_{\pa(i)\cap A}\!=\!a_{\pa(i)\cap A} )\label{eq:drop-predecessors-from-x-phil} \\
&= p(W_i \mid W_{\pa(i)\setminus A}, F_{\pa(i)\cap A}\!=\!a_{\pa(i)\cap A} ). \label{eq:drop-parents-from-x-phil} 
\end{align}
\end{lemma}

\noindent{\it Proof:} 
Here (\ref{eq:drop-successors-from-a-phil}) follows from Lemma \ref{lem:remove-b-from-joint-not-a-fn-phil} since by Definition
\ref{def:local-swig-mp-phil}, 
\[
W_{\pre(i)\cup \{i\}} \;\indep\; F_{A\setminus \pre(i)} \mid F_{A\cap \pre(i)}, F_A\neq \emptyset.
\]
Similarly, (\ref{eq:drop-predecessors-from-a-phil}) follows from Lemma \ref{lem:remove-b-from-conditional-not-a-fn-phil} since by the local Markov property:
\[
W_{i} \;\indep\; F_{(A\cap \pre(i))\setminus \pa(i)} \mid F_{A\cap \pa(i)}, F_{A\cap \pre(i)}\neq \emptyset.
\]
Finally, (\ref{eq:drop-predecessors-from-x-phil}) and (\ref{eq:drop-parents-from-x-phil}) again follow from the local Markov property
since 
\[
W_i \indep W_{\pre(i)\setminus \pa(i)}, W_{\pa(i)\cap A} \,\mid\,  W_{\pa(i)\setminus A}, F_A\!=\!a,
\]
hence $p(W_i \mid W_{\pre(i)}\!=\!w, F_A\!=\!a )$  does not depend on $w_{\pre(i)\setminus (\pa(i) \setminus A)} = 
(w_{\pre(i)\setminus \pa(i)}, w_{\pa(i)\cap A})$.  \hfill$\Box$\par

\subsection{Markov property for the observed distribution}

The following result shows that the Reformulated Local Markov property implies, via distributional consistency, the ordinary local Markov property for the observed distribution.
This result corresponds to Theorem \ref{thm:swig-obs-markov}.

\begin{theorem}
If the kernel $p(W_V \mid F_A)$ obeys distribution consistency and the augmented DAG local Markov property w.r.t. $\G^*$ then
$p(W_V)$ obeys the usual local Markov property w.r.t. $\G$.
\end{theorem}

\noindent{\it Proof:}  Let $w^* \in \mathfrak{X}_{\pre(i)}$.
\begin{align}
\MoveEqLeft{p(W_i=w \mid W_{\pre(i)}=w^*)}\nonumber\\
&= p(W_i=w \mid W_{\pre(i)}\!=\!w^*, F_{\pre(i)\cap A }=w^*_{\pre(i)\cap A})\nonumber\\
&= p(W_i=w \mid W_{\pa(i)\setminus A}\!=\!w^*_{\pa(i)\setminus A}, F_{\pa(i)\cap A }=w^*_{\pa(i)\cap A}).
\end{align}
Here the first equality follows by distributional consistency.
The second follows directly from the equality of (\ref{eq:drop-successors-from-a-phil}) and (\ref{eq:drop-parents-from-x-phil}) in Lemma \ref{lem:kernel-of-swig-lmp-phil}.
Since the last line only depends on $w^*_{\pa(i)}$, the ordered local Markov property for the DAG holds.
\hfill$\Box$

\subsection{Identifiability}

The next result shows that the Reformulated Local Markov property implies that the kernel $p(V \mid F_A)$ will be identified from the distribution of the observables
provided that the relevant conditional distributions are identified (from the distribution of the observables).
This result corresponds to Theorem \ref{thm:swig-identify}.

\begin{theorem}\label{thm:augmented-identify}
Suppose the kernel $p(W_V \mid F_A)$ obeys distribution consistency and the augmented DAG local Markov property w.r.t. $\G^*$.
Let $a$ be an assignment to the intervention targets in $A$, and let $v$ be an assignment to $W_V$. Then for every $i$:
\begin{align}
\MoveEqLeft{p(W_i=v_i \mid F_{A}=a, W_{\pre(i)} = v_{\pre(i)})}\nonumber\\
&= p(W_i=v_i \mid W_{\pa(i)\cap A}=a_{\pa(i)\cap A}, W_{\pa(i) \setminus A} = v_{\pa(i) \setminus A}).\label{eq:modularity}
\end{align}
Consequently, $p(W_V \mid F_A = a)$ is identified given $p(W_V)$ and obeys the Markov property for the DAG formed from $\G^*$ by removing all outgoing edges from vertices in $A$.
\end{theorem}

As before, we note that the equality (\ref{eq:modularity}) corresponds to the property referred to as `modularity' in the SWIG formulation which is also an instance
of the extended g-formula of \citet{robins86new,robins:effects:2004}.
\medskip

\noindent{\it Proof:} 
Let $a \in \mathfrak{X}_A$, $v\in \mathfrak{X}_V$. Now:
\begin{align}
\MoveEqLeft{p(W_i\!=\! v_i \mid W_{\pre(i)}\!=\! v_{\pre(i)}, F_A=a)}\nonumber\\[2pt]
&=  p(W_i\!=\! v_i \mid W_{\pa(i)\cap A}\!=\! v_{\pa(i)\cap A}, W_{\pa(i)\setminus A}\!=\! v_{\pa(i)\setminus A},\nonumber\\
&\kern180pt F_{A\cap \pa(i)}=a_{\pa(i)\cap A} ) \nonumber\\[2pt]
&=  p(W_i\!=\! v_i \mid W_{\pa(i)\cap A}\!=\! a_{\pa(i)\cap A}, W_{\pa(i)\setminus A}\!=\! v_{\pa(i)\setminus A},\nonumber \\
&\kern180pt F_{A\cap \pa(i)}=a_{\pa(i)\cap A} ) \nonumber\\[2pt]
&=  p(W_i\!=\! v_i \mid W_{\pa(i)\cap A}\!=\! a_{\pa(i)\cap A}, W_{\pa(i)\setminus A}\!=\! v_{\pa(i)\setminus A} ).
\end{align}
Here the first equality follows from the equality of (\ref{eq:swig-lmp-lhs-phil}) and (\ref{eq:drop-predecessors-from-x-phil});  the second follows from
the equality of (\ref{eq:drop-predecessors-from-x-phil}) and (\ref{eq:drop-parents-from-x-phil}); the third by distributional consistency.
\hfill$\Box$

\subsection{Distributions resulting from fewer interventions}

As in the SWIG case, a similar argument applies if we consider interventions on a subset $B \subseteq A$.
This result corresponds to Theorem \ref{thm:swig-gb}.

\begin{theorem}
Suppose the kernel $p(W_V \mid F_A)$ obeys distribution consistency and the augmented DAG local Markov property w.r.t. $\G^*$.
Let $b$ be an assignment to the intervention targets in $B\subseteq A$, and let $w^*$ be an assignment to $W_V$. Then for every $i$:
\begin{align}
\MoveEqLeft{p(W_i=w_i^* \mid F_{B}=b, W_{\pre(i)} = w^*_{\pre(i)}) 
}
\nonumber\\
&= p(W_i=w_i^* \mid W_{\pa(i)\cap B}=b_{\pa(i)\cap B}, W_{\pa(i) \setminus B} = w^*_{\pa(i) \setminus B}).\label{eq:modularity-for-b}
\end{align}
Consequently, $p(W_V \mid F_B = b)$ is identified given $p(W_V)$ and obeys the Markov property for the augmented DAG $\G^{**}$ formed from $\G^*$ by removing all outgoing edges from vertices in $B$ and removing the regime indicators $F_{A \setminus B}$.
\end{theorem}

\noindent{\it Proof:} 
\begin{align}
\MoveEqLeft{p(W_i =w_i \mid W_{\pre(i)}\!=\! w_{\pre(i)}, F_{B}\!=\!b) }\nonumber\\
&= p(W_i =w_i \mid W_{\pre(i)}\!=\! w_{\pre(i)}, F_{\pre(i)\cap B}\!=\!b_{\pre(i)\cap B})\nonumber\\[3pt]
&= p(W_i =w_i \mid W_{\pre(i)}\!=\! w_{\pre(i)}, F_{\pre(i)\cap B}\!=\!b_{\pre(i)\cap B},\nonumber\\[-2pt]
&\kern158pt F_{\pre(i)\cap (A\setminus B)}\!=\!w_{\pre(i)\cap (A\setminus B)})\nonumber\\[3pt]
&= p(W_i =w_i \mid W_{\pa(i)\setminus A}\!=\! w_{\pa(i)\setminus A}, W_{\pa(i)\cap B}\!=\! b_{\pa(i)\cap B},\nonumber\\[-2pt]
&\kern175pt W_{\pa(i)\cap ( A\setminus B) }\!=\! w_{\pa(i)\cap (A\setminus B)}) \nonumber\\[3pt]
&= p(W_i =w_i \mid W_{\pa(i)\setminus B}\!=\! w_{\pa(i)\setminus B},
W_{\pa(i)\cap B}\!=\! b_{\pa(i)\cap B}) \nonumber
\end{align}
Here the first equality is by Lemma \ref{lem:remove-b-from-joint-not-a-fn-phil}; the second is distributional consistency; the third follows from Theorem \ref{thm:augmented-identify} applied to $\G^*$; the fourth is a simplification.\hfill$\Box$

\section{The role of `fictitious' independence in Dawid's development}

\citet{dawid:2021}
 uses what he terms a `fictitious' independence
 in his proofs that the distribution of the kernels that condition on the regime indicators $F_i$ obey the Markov property for the augmented DAG with ITT variables. Specifically, in his proof of Lemma 4, though not the statement, he makes the formal assumption that 
\begin{align}
F_1 \indep F_0,\label{eq:fict1}
\end{align}
and similarly in the proof of Theorem 1 he assumes that all the regime indicators are mutually independent \citep[p.76, eqn. (82)]{dawid:2021}
\begin{align}
F_1 \indep F_2 \indep \cdots \indep F_{k-1} \indep F_k.\label{eq:fict2}
\end{align}
Such an independence assumption does not fit into the extended conditional independence (ECI) framework used by Dawid to describe the Markov property for augmented graphs. This is because, as stated by Dawid, an ECI statement $A \indep B \mid C$ must satisfy: ``(a) no non-stochastic variable occurs in $A$, and (b) all non-stochastic variables are included in $B\cup C$'' \citep[fn.~3]{dawid:2021}; these conditions allow independences to be viewed as well-defined restrictions on $p(A \mid B, C)$ since all of the non-stochastic variables appear on the right of the conditioning bar. However, an independence of the form $F_i \indep F_j$ violates both of these conditions.

Perhaps for this reason, Dawid argues that although his proofs make use of the assumptions (\ref{eq:fict1}) and (\ref{eq:fict2}) there is no loss of generality:
\begin{quotation}
So long as all our assumptions and conclusions are in the form described in footnote 3 [i.e. satisfy (a) and (b)], any proof that uses this extended understanding only internally will remain valid [\ldots] \citep[p.63]{dawid:2021}\par [\ldots]
because the premisses and conclusions of the argument relate only to distributions conditioned on the regime indicators, the extra assumption of variation independence is itself inessential, and can be regarded as just another ``trick.'' \citep[p.63, fn.24]{dawid:2021}
\end{quotation}

We will show via an example that Dawid's inference here is not valid: in general, the conclusion will not hold for a kernel without additional assumptions regarding the set of states taken by the non-stochastic variables. However,
notwithstanding this, as we also show below, Dawid's conclusions are still correct owing to the special structure that is present in the possible states taken by regime indicators.

\subsection{Invalid Implication}

To illustrate the issue with the proof we re-write Dawid's equations so as to make the argument transparent. 
Dawid makes the following claim:
\begin{claim}\label{claim:dawid}
Consider a kernel $q(x,y | a,b)$, with stochastic variables $X,Y$ and non-stochastic variables $A, B$. If the following ECI restrictions hold
\begin{align}
Y &\indep A \mid B,X\label{eq:1}\\
Y &\indep B \mid A,X\label{eq:2}
\end{align}
then it follows that:
\begin{align}
Y & \indep A,B \mid X.\label{eq:4}
\end{align}
\end{claim}

To relate this to Dawid's proof of Lemma 4, $A=F_0$, $B=F_1$, $X=X_0$ and $Y=\{H,Z,X_1^*\}$. Thus (\ref{eq:fict1}), (\ref{eq:1}), (\ref{eq:2}), (\ref{eq:4})
correspond to Dawid's Equations (41), (49), (50), (51), respectively.

In the proof of this claim, Dawid makes use of the `fictitious' independence $A \indep B$, but as noted above, he argues that this `internal' assumption may be made without loss of generality. To see that this implication does not hold without additional conditions on the state spaces for $A$ and $B$, suppose that the non-stochastic pair $(A,B) \in \mathfrak{S} \equiv \{-2,-1\}^2 \cup \{1,2\}^2$, so that $(A,B)$ take one of the following eight states:
\[
(-2,-2), (-2,-1), (-1,-2), (-1,-1), (1,1), (1,2), (2,1), (2,2).
\]
Note that, by construction, the non-stochastic variables are not variation independent; they always share the same sign.
Now let $P^{-}(Y,X)$ and $P^{+}(Y,X)$ be any pair of distributions over $(X,Y)$ such that $P^+(Y\mid X) \neq P^-(Y\mid X)$ and define the kernel
$p(Y, X \mid a,b)$ for $(a,b) \in \frak{S}$ as follows:
\begin{align}
p(Y, X \mid a,b) &= \left\{\begin{array}{cc}
p^{-}(Y,X) & \hbox{ if both } a,b < 0,\\
p^{+}(Y,X) & \hbox{ if both } a,b > 0,\\
\hbox{undefined} & \hbox{ otherwise}.
\end{array}
\right.
\end{align}
By construction, if $B=b<0$ then for all $a$ such that $(a,b) \in \mathfrak{S}$, i.e.~$a\in \{-2,-1\}$ it holds that
\begin{align}
p(Y \mid A=a, B=b, X) = p^-(Y\mid X),
\end{align}
so (\ref{eq:1}) holds when $B=b<0$. The argument when $B=b >0$ is symmetric, since in this case for all $a$ such that 
 $(a,b) \in \mathfrak{S}$, we have 
$p(Y \mid A=a,B=b, X) = p^+(Y\mid X)$.
Hence (\ref{eq:1}) holds for all $B \in \{-2,-1,1,2\}$. A symmetric argument replacing $A$ with $B$ shows that (\ref{eq:2}) also holds.

However, the conclusion (\ref{eq:4}) fails since by construction:
\begin{align}
p(Y \,|\, A\!=\!-2,B\!=\!-2, X) &= p^-(Y\,|\, X)\nonumber \\
 &\neq p^+(Y\,|\, X) \;=\; p(Y \,|\, A\!=\!2,B\!=\!2, X).
\end{align}
 
The implication in the claim corresponds to an ECI instance of the Intersection Axiom (CI5) introduced by Dawid in his classic papers \citep{dawid:1979,dawid:1980}. As he notes in several places in \citep{dawid:2021}, this implication is well known not to hold in general.%
\footnote{For recent related work giving general conditions under which this implication holds for ordinary (not extended) conditional independence, 
see \citet{gill:2019}, \citet[Ch.4]{sullivant:2018}, \citet{peters:intersection}.} Our counterexample above simply serves to show that even though there is no distribution over the non-stochastic variables, the implication will not hold if the non-stochastic variables are not variation independent.

\subsection{Validity of the Conclusion for Regime Indicators}

That the implications used in Dawid's proofs of Lemma 4 and Theorem 1 do not hold  -- without conditions on the joint space for the regime indicators --  may at first seem to call into question Dawid's conclusions.  However, at least in causal theories making use of DAG representations and involving multiple treatments, the decisions as to whether to intervene, and if so, which value to enforce are unconstrained. Consequently variation independence will hold, and hence the conclusion will be valid.

However, there are situations in which interventions may be constrained. For example suppose that there are two strategies for a medical condition;  each treatment involves two separate stages $(A_1,A_2)$. At time $t=1$, the doctor must decide between strategies `$1$', `$2$' .
It is easy to imagine situations in which, if treatment was commenced at time $1$,  the treatment at time $2$ involves `completing' the treatment that was started at time $1$; for example removing surgical stitches from the specific operation performed at time $1$.
In this case, the treatment options available at time $2$ are constrained by the decision at time $1$.

Reflecting this, there have been causal decision theories proposed in which variables do not live in a product space; see \citet{THWAITES2010889}. 
Likewise, in the potential outcome framework, the formulation of Causally Interpreted Structured Tree Graphs given by \citet{robins86new} also allows for this possibility.

However, even in this case, Dawid's implication will still hold, provided that the following condition holds:

\begin{definition}
Let $\mathfrak{F}_A \subsetneq \bigtimes_{i \in A} \mathfrak{X}_{i}\cup\{\emptyset\}$ indicate the (constrained) state-space for the set of regime indicators $F_A$. 
\begin{align}
\MoveEqLeft{\hbox{For all } f \in \mathfrak{F}_A, \hbox{ s.t.~}f\neq \emptyset\;}\\
& \Rightarrow\; \nonumber
\hbox{\it there exists } i\in A \hbox{ \it s.t.~} f_i \neq \emptyset \hbox{ \it and } (f_{-i}, \emptyset) \in \mathfrak{F}_A,\label{eq:move-to-idle}
\end{align}
where $f_{-i}$ indicates the values assigned to $A\setminus \{i\}$ by $f$.
\end{definition}
In words, this states that for any possible setting of the regime indicators, in which they are not all `idle', there exists some intervention target $A_i$ that is intervened upon under $f$, that could have not been intervened upon, such that the resulting vector $(f_{-i}, \emptyset)$ is still a valid value for $F_A$.

This condition may still hold in settings in which, if a later target is intervened upon, the regime under which an earlier target is set to `idle' is not well-defined. For example, an intervention on $A_2$ setting $F_2=1$ may only be well-defined if $F_1=1$, but not $F_1=\emptyset$.
In the treatment completion example above this would be the case 
if, in the absence of an intervention on $A_1$, some patients would receive treatment $2$ at time $1$, so that the subsequent intervention $F_2=1$ would not be well-defined. If the same holds for $F_2=2$ then, $F_1=\emptyset \Rightarrow F_2=\emptyset$.\footnote{Here we are implicitly supposing that only static regimes are under consideration so that $F_2$ can only take the values $1$, $2$ or $\emptyset$.}
 The condition (\ref{eq:move-to-idle}) will always hold provided that treatment decisions follow a time order, and that, regardless of the decisions that have occurred previously, it is always possible to decide to replace the `last' intervention with the idle regime.

It is easy to see that under the condition (\ref{eq:move-to-idle}) for any $f\in  \mathfrak{F}_A$, there will exist a sequence $(f=f^0,f^1,\ldots ,f^q=\emptyset)$ such that for $j=1,\ldots , q$, $f^j \in \mathfrak{F}_A$, and $f^{j}$ contains one more idle regime indicator than $f^{j-1}$. It then follows under this condition that:
\begin{align}
\hbox{if for all } i\in A,\quad W \indep F_i \mid F_{A\setminus \{i\}}\quad \hbox{ then } \quad W \indep  F_{A},
\end{align}
where here the conditional independence statements implicitly quantify over all the assignments to $F_A$ that are in $\mathfrak{F}_A$ and hence valid.

\eject

\subsubsection*{Funding Information} The authors completed work on this paper while visiting
the American Institute for Mathematics and the Simons Institute, Berkeley.
The authors were supported in by the ONR grant N000141912446; Robins was also supported by
NIH grant R01 AI032475.

\subsubsection*{Conflicts of Interest}  The authors state no conflicts of interest.

\subsubsection*{Acknowledgments} We thank Ilya Shpitser and Philip Dawid for helpful comments and discussions.
 
 \bibliographystyle{chicago}
\bibliography{dawid-discussion}

\appendix

\section{Appendix}

\subsection{Conditions implying the SWIG local Markov property for $\G$ given that $p(V)$ factors with respect to $\G$
 }\label{app:dawid-swig}

Here we show that if ${\cal P}_A$ obeys the SWIG local Markov property corresponding to a complete graph $\overline{\G}$,  and further,  the observed distribution   
$p(V)$ is positive and obeys the local Markov property for a subgraph $\G$ of $\overline{\G}$ then it follows from distributional consistency that ${\cal P}_A$ also obeys the SWIG local Markov property corresponding to ${\G}$.

\begin{theorem}\label{thm:top-and-bottom-swig}
Suppose ${\cal P}^{\subseteq}_A$ obeys distributional consistency and ${\cal P}_A$  obeys the 
SWIG ordered local Markov property for $\overline{\G}$, a complete DAG.\footnote{A DAG is complete if there is an edge between every pair of variables. Note that in this case
there is only one topological ordering.} 
If $p(V)$ is positive and obeys the Markov property for a subgraph $\G$ of $\overline{\G}$, then 
${\cal P}_A$ obeys the SWIG ordered local Markov property for $\G$.
\end{theorem}

\noindent{\it Proof:} Let $v \in \mathfrak{X}_{\pre(i)}$, $v_i^* \in \mathfrak{X}_{i}$ and $a \in \mathfrak{X}_A$.
\begin{align}
\MoveEqLeft{p(X_i({ a})=v^*_i \mid X_{\pre(i)}({ a})= v_{\pre(i)})}\label{eq:first-bolt-on-eq}\\
&= p(X_i=v^*_i \mid  X_{\pre(i) \setminus A} = v_{\pre(i) \setminus A}, X_{\pre(i)\cap A}=a_{\pre(i)\cap A})\label{eq:first-consistency-AA}\\
&= p(X_i=v^*_i \mid  X_{\pa(i) \setminus A} = v_{\pa(i) \setminus A}, X_{\pa(i)\cap A}=a_{\pa(i)\cap A}).\label{eq:first-consistency-BB}
\end{align}

Here the first equality follows from (\ref{eq:modularity-swig}) which holds under the local Markov property for $\overline{\cal G}$.
The second equality is due to the local Markov property for $p(V)$. Consequently, we see that (\ref{eq:first-bolt-on-eq})
does not depend on $v_{\pre(i) \setminus (\pa(i) \setminus A)}$ nor on  $a_{A \setminus \pa(i)}$
 as required by the SWIG local Markov property. Note that positivity is used here in order to ensure that (\ref{eq:first-consistency-AA}) and
 (\ref{eq:first-consistency-BB}) are equal for all assignments to the variables in the conditioning events.
\hfill$\Box$
\medskip

This result is similar in spirit to Dawid's construction in that it provides conditions that, in conjunction with the observed distribution $p(V)$ obeying the Markov property for $\G$, are sufficient to imply that ${\cal P}_A$ obeys the SWIG local Markov property for $\G$. The SWIG ordered local Markov property on ${\cal P}_A$ for a complete graph $\overline{\G}$ corresponds to the Finest Fully Randomized Causally Interpreted Structured Tree Graph (FFRCISTG) of \citet{robins86new}, in the case where $A$ represents the finest, i.e. largest, set of treatment variables for which well-defined counterfactuals exist, and there are no (population or individual level) exclusion restrictions.

Note that for a complete graph $\overline{\cal G}$, for every variable $i$, $\pre(i) = \pa(i)$. Consequently, 
the SWIG local Markov property (\ref{eq:d-sep-local-swig2}) and (\ref{eq:d-sep-local-swig4}) reduce to requiring that for every $i$: 
\begin{align}\label{eq:d-sep-local-swig2-comp}
 X_i({a}) \indepd 
 \overbrace{X_{\pre(i) \cap A}(a)}^{\hbox{\tiny ignorability}},\;\;
 \overbrace{\vphantom{X}x_{A\setminus \pre(i)}}^{\hbox{\tiny time order}}\;\;   
\mid 
\underbrace{x_{A \cap \pre(i)},}_{\hbox{\tiny~~ \parbox{40pt}{\centering fixed predecessors}~~}}
\underbrace{X_{\pre(i) \setminus A}(a).}_{\hbox{\tiny ~~\parbox{40pt}{\centering random predecessors}~~}}
 \end{align}
 Thus we see that the SWIG local Markov property for  the complete graph $\overline{\cal G}$ solely imposes ignorability and that interventions in the future do not change (the distribution of) variables in the past. 

\medskip

The single graph approach given by Definition \ref{def:local-swig-mp}
and the two graph construction of Theorem \ref{thm:top-and-bottom-swig} each have their own strengths and weaknesses:

\begin{itemize}
\item In the single graph approach the model places restrictions on ${\cal P}_A$; distributional consistency for
${\cal P}^{\subseteq}_A$ then implies the relevant SWIG Markov properties for all the other distributions in ${\cal P}^{\subseteq}_A$, including the factual distribution $p(V)$. This approach is more concise insofar as it requires fewer conditions. The approach does not require $p(V)$ to be positive.
\item In the two graph construction, the graph $\G$ specifies conditional independence restrictions on the observed distribution $p(V)$ via an ordinary Markov property, while the SWIG Markov property for the complete supergraph $\overline{\cal G}$ imposes ignorability and a total time order on ${\cal P}_A$. Under positivity for $p(V)$, distributional consistency for ${\cal P}^{\subseteq}_A$ then implies the relevant SWIG Markov properties for every distribution in ${\cal P}^{\subseteq}_A$. Though it requires more conditions this approach has the advantage that it clearly demarcates a set of additional  conditions that, when added to the assumption that $p(V)$ obeys the Markov property for $\G$, suffice to construct the full model on ${\cal P}^{\subseteq}_A$.
\end{itemize}

The fact that the single graph approach does not require positivity can be seen as an advantage since
it does not restrict the set of observed distributions. As a consequence, the graph in the single graph approach may include edges that indicate effects arising from interventions on $A$ that set variables to configurations that have probability zero under the observed distribution. Even in the absence of confounding, such effects may only be detectable via randomized experiments; see \cite{robins:mediation:2020} for further discussion.

\subsection{Derivation of part of the augmented DAG local Markov property for $\G$ from $p(V)$}

Similar to our development in \S \ref{app:dawid-swig}, and also to Dawid's construction, we provide conditions on the kernel $p(W_V \mid F_A)$ that, in conjunction with 
a positive observed distribution distribution  
$p(V)$ that obeys the local Markov property for a subgraph $\G$, suffice to ensure $p(W_V \mid F_A)$
obeys the local property for the corresponding augmented graph ${\G}^*$. These conditions are formulated in terms of a decision diagram $\overline{\G}^*$
corresponding to a complete DAG  $\overline{\G}$ that contains $\G$ as a subgraph.

\begin{theorem}
Suppose that $p(W_V \mid F_A)$ obeys distribution consistency and the augmented DAG local Markov property with respect to $\overline{\G}^*$ where $\overline{\G}$ is a complete DAG. If $p(V)$ is positive and obeys the (ordinary) Markov property for a subgraph $\G$ of $\overline{\G}$ then 
$p(W_V \mid F_A)$ also obeys the augmented DAG local Markov property for ${\G}^*$.
\end{theorem}

For a complete graph $\overline{\cal G}$, for every variable $i$, $\pre(i) = \pa(i)$, and thus the
 local Markov property for $\overline{\cal G}^*$ requires that for every $i$,
\begin{align}\label{eq:d-sep-local-swig2-comp-aug}
W_i  \indepd 
 \overbrace{W_{\pre(i) \cap A}}^{\hbox{\tiny ignorability}},\;\;
 \overbrace{\vphantom{W}F_{A\setminus \pre(i)}}^{\hbox{\tiny time order}}\;\;   
\mid 
\underbrace{F_{A \cap \pre(i)},}_{\hbox{\tiny~~ \parbox{40pt}{\centering fixed predecessors}~~}}
\underbrace{W_{\pre(i) \setminus A},}_{\hbox{\tiny ~~\parbox{40pt}{\centering random predecessors}~~}}
\underbrace{F_{A\vphantom{\mid}} \neq \emptyset.}_{\hbox{\tiny \parbox{40pt}{\centering intervene on all of $A$}}}
 \end{align}
 Thus, similar to (\ref{eq:d-sep-local-swig2-comp}), this imposes ignorability and that interventions in the future do not change (the distribution of) variables in the past. 
 \bigskip

\noindent{\it Proof:} Let $v \in \mathfrak{X}_{\pre(i)}$, $v^* \in \mathfrak{X}_{i}$ and $a \in \mathfrak{X}_A$.

\begin{align}
\MoveEqLeft{p(W_i=v^*_i \mid F_{A}=a, W_{\pre(i)} = v_{\pre(i)})}\label{eq:first-bolt-on-eq-aug}\\
&= p(W_i=v^*_i \mid W_{\pre(i)\cap A}=a_{\pre(i)\cap A}, W_{\pre(i) \setminus A} = v_{\pre(i) \setminus A})\label{eq:first-consistency-AA-aug}\\
&= p(W_i=v^*_i \mid W_{\pa(i)\cap A}=a_{\pa(i)\cap A}, W_{\pa(i) \setminus A} = v_{\pa(i) \setminus A}).\label{eq:first-consistency-BB-aug}
\end{align}

Here the first equality follows from (\ref{eq:modularity}) which holds under the augmented local Markov property for $\overline{\cal G}^*$.
The second equality is due to $p(V)$ obeying the local Markov property for $\G$. Consequently, we see that (\ref{eq:first-bolt-on-eq-aug})
does not depend on $v_{\pre(i) \setminus (\pa(i) \setminus A)}$ nor on  $a_{A \setminus \pa(i)}$,
 as required by the augmented graph local Markov property. Note that positivity ensures that (\ref{eq:first-consistency-AA-aug}) and
 (\ref{eq:first-consistency-BB-aug}) are equal for all assignments to the variables in the conditioning events.
\hfill$\Box$

\end{document}